\documentclass[a4paper,twoside,english]{scrartcl}
\usepackage[T1]{fontenc}
\usepackage[latin1]{inputenc}
\usepackage{babel}
\usepackage{verbatim}
\usepackage{amsthm}
\usepackage{amsmath}
\usepackage{amssymb}
\usepackage{scrpage2}
\usepackage[unicode=true,pdfusetitle,
 bookmarks=true,bookmarksnumbered=false,bookmarksopen=false,
 breaklinks=false,pdfborder={0 0 1},backref=false,colorlinks=false]
 {hyperref}
\usepackage{breakurl}


\usepackage{a4wide}
\usepackage{tu-preprint}

\usepackage{amssymb}
\usepackage{amsmath}
\usepackage{amsthm}
\theoremstyle{plain}
\usepackage{graphicx}
\usepackage{epstopdf}
\usepackage{topcapt}
\usepackage{latexsym}
\usepackage{graphics}
\usepackage{subfigure}
\newenvironment{keywords}{ \noindent\footnotesize\textbf{Keywords and phrases:}}{}

\newenvironment{class}{\noindent\footnotesize\textbf{Mathematics subject classification 2010:}}{}

\numberwithin{equation}{section}

\newtheorem{theorem}{Theorem}[section]

\newtheorem{corol}[theorem]{Corollary}
\newtheorem{prop}[theorem]{Proposition}

\newtheorem{remark}[theorem]{Remark}

\newtheorem{lemma}[theorem]{Lemma}
\newtheorem{definition}[theorem]{Definition}

\def\b{\beta}

\def\cF{{\cal F}}

\def\d{\partial}
\def\de{\delta}
\def\e{\varepsilon}

\def\8{\infty}

\def\N{\mathbb N}

\def\och{{\overline \chi}}

\def\Ome{\Omega}
\def\ops{{\overline \psi}}
\def\os{{\overline s}}

\def\ot{{\overline \tau}}

\def\oB{\overline {B(x_2,b)}}

\def\obeta_2{{\overline \beta_2}}

\def\oy{{\overline y}}
\def\oz{{\overline z}}

\def\0*{^{\odot *}}

\def\R{\mathbb{R}}

\def\s{\sigma}

\def\t{\tau}
\def\th{\theta}

\def\vi{\varphi}

\def\Proof{\noindent{\em Proof. }}
\def\<{\langle}
\def\>{\rangle}
\def\_>{\longrightarrow}
\begin{document}

\institut{Institut f\"ur Analysis}

\preprintnumber{MATH-AN-07-2014}

\preprinttitle{A differential equation with state-dependent delay
from cell population biology}

\author{Philipp Getto, Marcus Waurick}

\makepreprinttitlepage\setcounter{section}{-1}

\date{}

\title{A differential equation with state-dependent delay
from cell population biology}

\author{Philipp Getto\\ Philipp.Getto@tu-dresden.de\\ Marcus Waurick\\ Marcus.Waurick@tu-dresden.de \thanks{BCAM (Basque Center for Applied Mathematics), 48009 Bilbao, Spain, and TU Dresden, Fachrichtung Mathematik, Institut f\"ur Analysis, 01062 Dresden, Germany, Tel.: ++49 351 463 34054 }}
\maketitle
\begin{abstract}
\noindent \textbf{Abstract.} We analyze a differential equation with a state-dependent delay that is implicitly
defined via the solution of an ODE. The equation describes an established though little 
analyzed cell population model. Based on theoretical results of Hartung, Krisztin, Walther 
and Wu we elaborate conditions for the model ingredients, in particular vital 
rates, that guarantee the existence of a local semiflow. Here proofs are based on implicit function arguments. 
To show global existence, we adapt a theorem from a classical book on functional differential equations by Hale and Lunel, which gives conditions under which - if there is no global existence - closed and bounded sets are left for good, to the $C^1$-topology, which is the natural setting when dealing with state-dependent delays. The proof is based on an
older result for semiflows on metric spaces.  
\end{abstract}

\begin{class}  34K05; 92C37; 37N25
\end{class}\setcounter{page}{1}

\newpage
\tableofcontents{}

\newpage

\section{Introduction}
\begin{figure}
\centering{}
\includegraphics[scale=0.46]{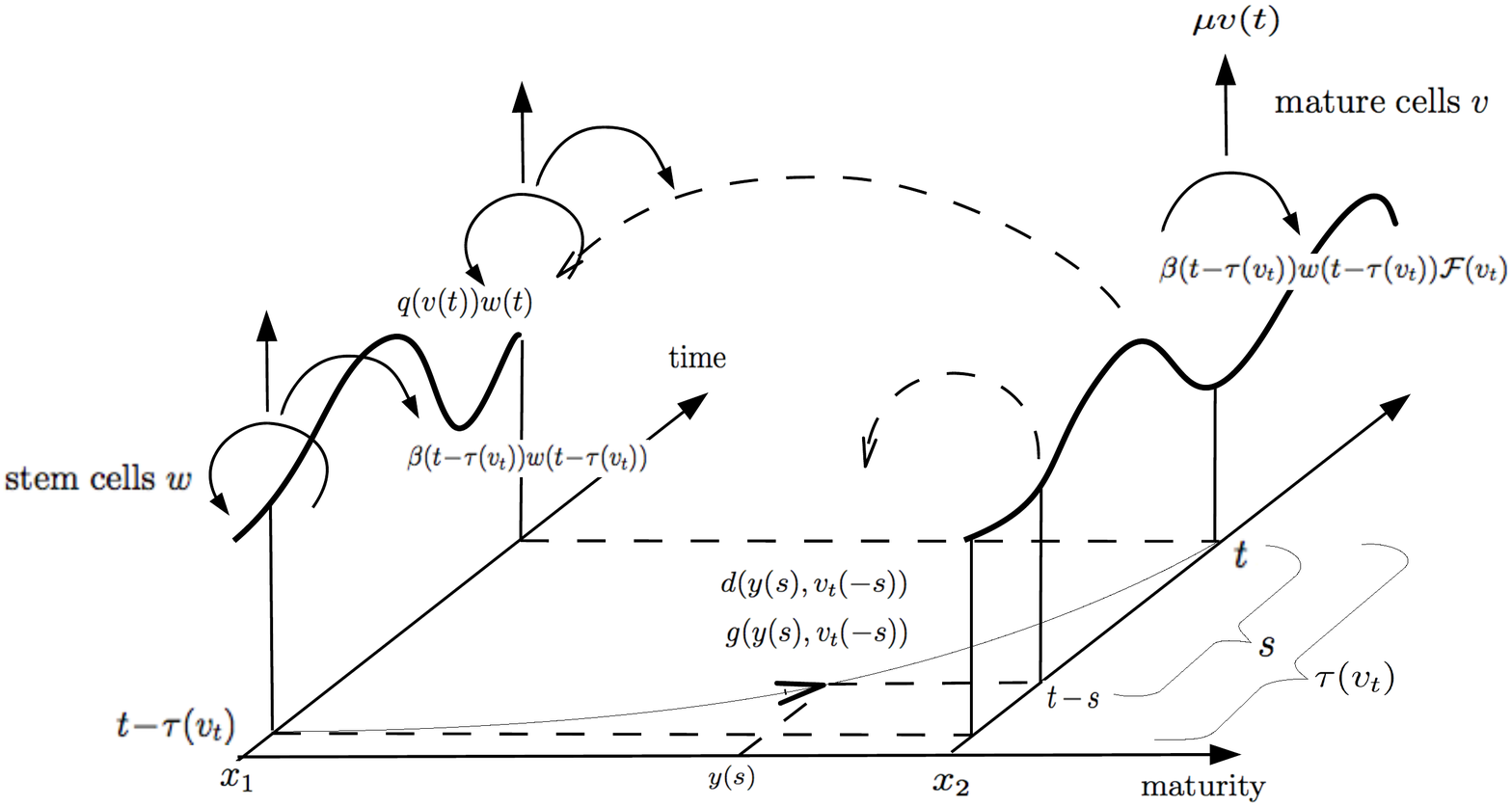}
\caption{
\label{f1}
Maturation process of stem cells modeled as differential equations with state-dependent
delay: Dashed arrows refer to regulation by mature cells. Continuous arrows refer to mortality (vertical), self-renewal (circular, anti-clockwise) and transition to higher maturity
(circular, clockwise).}
\end{figure}

In this paper we analyze a class of differential equations of the form
\begin{eqnarray}
w'(t)&=&q((v(t))w(t), 
\label{eq68}\\
v'(t)&=&\b(v(t-\t(v_t)))w(t-\t(v_t))\cF(v_t)-\mu v(t). 
\label{eq11}
\end{eqnarray}
\noindent
We use the standard notation
\[
x_t(s):=x(t+s),\;\;s<0,
\] 
if a function $x$ is defined in $t+s\in\mathbb{R}$. If $t$ is fixed, then $x_t$ is a function describing the history of $x$ at time $t$. Both (\ref{eq68}) and (\ref{eq11}) are equations in 
$\mathbb{R}$ and all functions are real-valued. Next, $\t$ and $\cF$  are nonlinear functionals with nonnegative values and domain in a space of functions, $\b$ and $q$ denote functions with real arguments and $\mu$ is a parameter. 
The functional $\t$ describes the delay and is allowed to depend exactly on the second component $v_t$ of the two components of the system at time $t$. The delay is in general only
implicitly given. We specify $\t$ in $v_t$ as the solution of the equation $y(\t,v_t)=x_1$, where
$y(v_t)=y(\cdot,v_t)$ is defined via the ODE
\begin{eqnarray}
y'(s)&=&-g(y(s),v_t(-s)),\;s>0,
\nonumber\\
y(0)&=&x_2,
\nonumber
\end{eqnarray}
where $x_1, x_2\in\mathbb{R}$, $x_1<x_2$ are given model parameters and $g$ is a 
given (nonnegative) model function, see Figure \ref{f1}. The functional $\cF$ can, in $v_t$, be specified as
\begin{eqnarray}
\cF(v_t):=e^{\int_0^{\t(v_t)}d(y(s,v_t),v_t(-s))ds},
\label{eq81}
\end{eqnarray}
where $d$ is another given (nonnegative) model function. 
Equations (\ref{eq68}--\ref{eq11}) together with the ODE can be classified as a differential equation with implicitly defined delay with state dependence. 

The system describes the maturation process of stem cells. 
It was first formulated in \cite{Alarcon}  following an idea of Anna Marciniak-Czochra to analyze a partial differential equation of transport type, that is a variant of ODE models for related biological problems. In the following we summarize the biological interpretation and refer to \cite{Alarcon, Doumic} and references given therein for further biological background and modeling aspects.
See Figure \ref{f1} for visualization of the exposition.

The dynamics of the whole cell population can be described in terms of the dynamics of the concentration of stem cells $w$ and fully mature cells $v$. The function $q$ is the stem cell population net growth rate. This rate describes stem cell population growth due to division as well as outflow due to maturation
 or decay.  It is regulated by the current size of the mature cell population. 
 The function $\b$ is a rate describing the outflow of those stem cells that commit themselves to maturation. It is also regulated. Next, $\t(v_t)$ is the time it takes until a stem cell committed to maturation becomes fully mature given that it survives and that at the moment of entering the mature cell compartment it is time $t$ and a history of regulation by mature cells $v_t$ is experienced. The value $\cF(v_t)$ can be interpreted as the {\it population net growth factor} during the maturation process, which means that if $\Psi$ is a given outflow of committed stem cells  at time $t-\t(v_t)$ then $\cF(v_t)\Psi$ is the inflow of cells into the mature cell compartment $\t(v_t)$ time units later, i.e., at time $t$. Finally $\mu$ is the decay rate of mature cells. This completes our interpretation of (\ref{eq68}--\ref{eq11}).  

In a further specification, we describe the maturity of a cell by a one-dimensional variable $x\in\R$. We assume that maturation occurs at a rate $g(x,y)$ that 
depends on maturity $x$ and on the current size of the population of mature cells $y$. The stem cells $w$ are then cells at  initial maturity $x=x_1$, 
and  the fully mature cells $v$ are cells of maturity $x=x_2>x_1$. We here do not describe 
the population size development for progenitor cells, i.e., cells with maturities 
$x\in(x_1,x_2)$, though, in general, it can be computed, see \cite{Alarcon}. Moreover, we now assume that the net growth of the cell population 
during maturation, that includes reproduction
and decay of cells, can be described by a per cell net rate $d=d(x,y)$ that, like $g$, 
depends on maturity and mature cell population.
This leads to the given specification of $\cF$. 
%

In \cite{Walther}, Hartung et al. elaborated conditions that can be satisfied by differential equations with state-dependent delay and that  guarantee the existence of a local (in time) semiflow for such equations. The main idea is to restrict initial histories to a submanifold of a space of $C^1$-functions. 
In \cite[Proposition VII 2.2]{3Diekmann} the authors present a criterion for global existence for local semiflows on metric spaces. The idea is to show that if the maximal time interval of existence is finite, then an arbitrary compact set at some point in time is left for good, see Figure \ref{f2} (left). 
In applications it is useful to have a variant of this criterion in which the assumption of compact sets is relaxed to closed and bounded sets, see Figure \ref{f2} (right). For non-autonomous functional differential equations with functionals that satisfy smoothness conditions in a setting of continuous functions, such a criterion is \cite[Theorem 2.3.2]{Hale}. 


In this paper we apply the mentioned results of \cite{Walther, 3Diekmann} to establish local and global existence for (\ref{eq68}--\ref{eq11}). As our main result we consider the elaboration of  respective sets of conditions on the rates $\b$, $q$, $g$ and $d$ that guarantee local and global existence of (\ref{eq68}--\ref{eq11}). As our aim is to preserve generality where possible, we do this stepwise in a top down approach. 

At the top level we adapt the criterion for global existence of \cite{3Diekmann} to the setting for state dependent delay equations of \cite{Walther} and obtain a new sufficient criterion for global existence, similar to the mentioned one in \cite{Hale}, but applicable to general differential equations with state-dependent delay. 

We then study a class of delay differential equations (DDE) of the form
\begin{eqnarray}
x'(t)=A(x_t)x(t)+b(x_t),
\label{eq69}
\end{eqnarray}
where $A$ is a diagonal matrix valued functional and $b$ a vector valued functional.
Note that this class contains (\ref{eq68}--\ref{eq11}). For (\ref{eq69}) we can show a variation of constants formula with which we establish useful bounds for the trajectory. We show local and global existence for this class of DDE. 

The next step is to elaborate conditions for $\b$, $q$, $\t$ and $\cF$ such that we can  
use the previously established theory for (\ref{eq69}). These are differentiability
and Lipschitz conditions in finite dimensions, in the case of $\b$ and $q$, and in infinite dimensions, in the case
of $\t$ and $\cF$. 

Finally we establish properties for $g$ and $d$ that guarantee that the conditions for $\t$ and $\cF$ hold. This amounts to defining, often implicitly, nonlinear operators and showing their differentiability with the implicit function theorem. 

In order to highlight our results before the technicalities and to come soon to our main theorem (Theorem \ref{theo5}), we have opted for the following structure for the remainder of the paper.  
Section \ref{s1} is devoted to the presentation of both existing results from the literature 
and our main results in a precise mathematical setting and in
Section \ref{s2} we elaborate proofs of our main results. In each of these two sections 
each subsection refers to one step of the discussed top down approach. In particular 
Section \ref{S1} contains the proofs of Section \ref{ss1},  Section \ref{ss5} the proofs of 
Section \ref{ss2}, etc..  
We close the paper with a discussion section. 
%
%
\section{Existing results and main results of the paper}\label{s1}
%
%
\subsection{Differential equations with state dependent delay}\label{ss1}
%
%
Conditions for existence and uniqueness of a noncontinuable  solution for differential
equations with state dependent delay are given in \cite{Walther}. We start by summarizing
these results. For $n\in\mathbb{N}$ we will use the Banach spaces
\begin{eqnarray}
(C([a,b],\mathbb{R}^n),\|\cdot\|),\;\|\phi\|:=\max_{\th\in[a,b]}|\phi(\th)|,\;
\nonumber\\
\;(C^1([a,b],\mathbb{R}^n),\|\cdot\|^1),\;\|\phi\|^1:=\|\phi\|+\|\phi'\|.
\nonumber
\end{eqnarray}
Next,  we define $C:=C([-h,0],\R^n)$ and $C^1:=C^1([-h,0],\R^n)$ for some $h\in(0,\infty)$. 
Let 
\[
U\subset C^1\;{\rm open},\;\;f:U\_>\mathbb{R}^n. 
\]
Then we can define solutions for DDE:
\begin{definition}\rm 
For any $\phi\in U$, a {\it solution} on $[-h,t_*)$, for some $t_*\in(0,\infty]$, of the initial value problem (IVP)
\begin{eqnarray}
x_0=\phi,\;\;x'(t)=f(x_t),\;\;t>0
\label{eq97}
\end{eqnarray}
is a continuously differentiable function $x:[-h,t_*)\longrightarrow\mathbb{R}^n$, which satisfies $x_t\in U$ for all $t\in (0,t_*)$ as well as the IVP.
\end{definition}
\noindent
Solutions on closed intervals $[-h,t_*]$, $t_*>0$ are defined analogously. 
A necessary condition for the unique solvability of IVPs is that initial data are restricted to the closed set 
\begin{eqnarray}
X=X(f):=\{\phi\in U:\phi'(0)=f(\phi)\}.
\label{eq96}
\end{eqnarray}
So naturally one requires that $f$ is chosen in a way that
$X=X(f)$ is nonempty. 
Moreover the following smoothness condition (S) is appropriate: 
\begin{itemize}
\item[(S1)] $f: U\_> \mathbb{R}^n$ is continuously differentiable,
\item[(S2)] each derivative $Df(\phi)$, $\phi\in U$ extends to a linear map 
$D_e f(\phi):C\_>\mathbb{R}^n$ and
\item[(S3)] the following map is continuous
\[
U\times C\_>\mathbb{R}^n,\;\;(\phi,\chi)\longmapsto (D_ef)(\phi)\chi.
\] 
\end{itemize}
We can then rephrase parts of Theorem 3.2.1 in \cite{Walther} as
\begin{theorem}\label{theo2} {\rm (Local semiflow)} Suppose that $f$ satisfies $(S)$ and is
such that $X$ is nonempty. Then $X$ is a continuously 
differentiable submanifold of $U$ with codimension $n$. Moreover, for each $\phi\in X$ 
there exists some $t_\phi>0$ and a unique noncontinuable solution $x^\phi:[-h,t_\phi)\_>\mathbb{R}^n$ of the IVP.
All segments $x_t^\phi$, $t\in [0,t_\phi)$, belong to $X$
and for
\[
\Omega:=\{(t,\phi):\;t\in [0,t_\phi),\; \phi\in X\}
\]
the map
\begin{eqnarray*}
S:\Omega\_> X;\;\;S(t,\phi):=x_t^\phi
\end{eqnarray*}
defines a continuous semiflow.
\end{theorem}
By the {\it existence of a local semiflow} we shall mean that the conditions in the conclusions of the previous theorem hold. For $\phi\in X$, we denote by $I_\phi:= [0,t_\phi)$  {\it maximal intervals} of existence of $x_t^\phi$.

\begin{figure}
\centering{}
\includegraphics[scale=0.2]{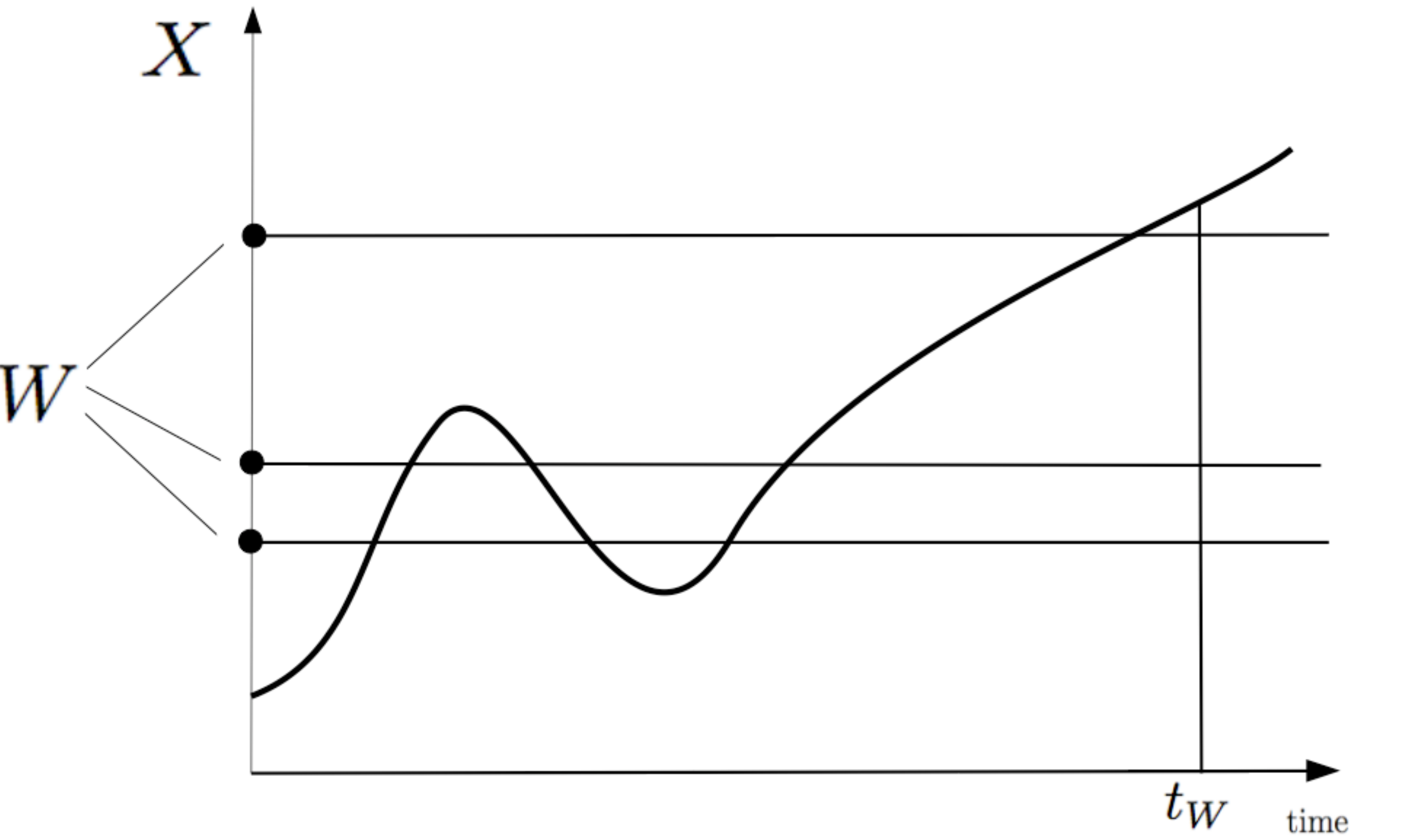}
\includegraphics[scale=0.2]{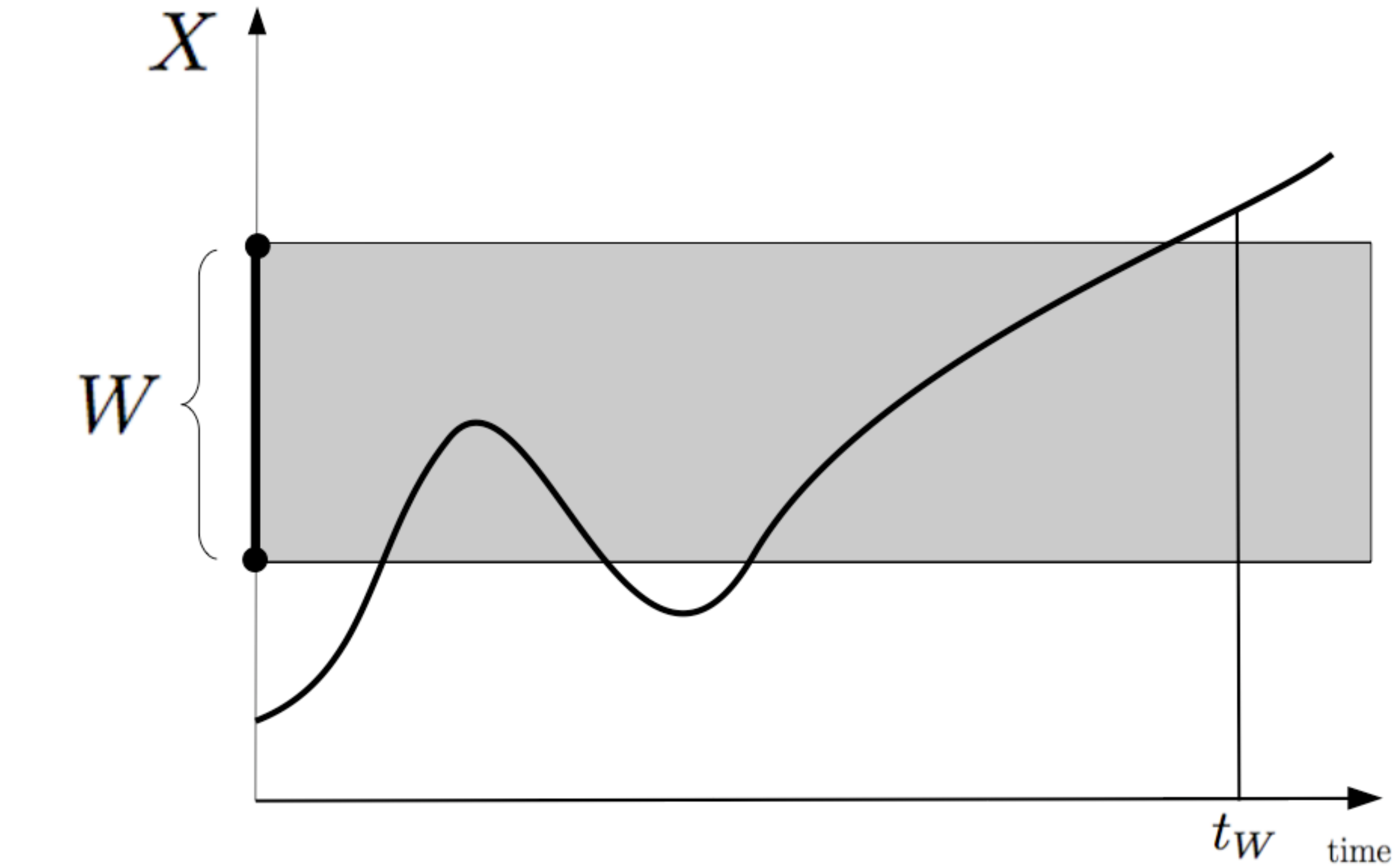}
\caption{
\label{f2}
Schematic visualization of sufficient criteria for global existence: an arbitrary compact set $W$ (left), which in infinite dimensions is a set with empty interior, and a closed and bounded set $W$ (right) are left for good at some point $t_W$.}
\end{figure}

 Proposition VII 2.2 in \cite{3Diekmann} states some properties of semiflows on 
metric spaces in a context in which completeness of the metric
space is assumed. See also the earlier Section II 10 in Amann \cite{Amann} on flows. The following result is an application of part (iii) of the proposition to $X$.
We do not know whether $X$ is complete (we know that it is closed in 
the relative topology of the open set $U$) but for Proposition VII 2.2 (iii) in \cite{3Diekmann} 
completeness is not necessary and the proof needs not to be changed if one drops the 
completeness assumption, see also Figure \ref{f2} (left). 
\begin{lemma}\label{lem3}
Suppose that there exists a local semiflow on $X$. 
Let $\phi\in X$ and assume that $t_\phi<\infty$. Then, for any $W\subset X$ compact there exists some $t_W$, such that $x_t\notin W$ 
for all $t\in[t_W,t_\phi)$. 
\end{lemma}
\noindent
For non-autonomous functional differential equations with functionals that are defined on $C$-open subsets of $C$ and  $C$-completely continuous  
Theorem 2.3.2 in \cite{Hale} is a variant of this statement in which the assumption of 
compact sets is relaxed to closed and bounded sets, see Figure \ref{f2} (right). Also in our setting we obtain such a result if we add more smoothness assumptions:
\begin{definition}\rm
A functional $A:U\_>\R^{n\times m}$ is called (Lb) (where L stands for Lipschitz and b for bounded) if for any $C^1$-bounded set $B\subset U$ there exists some $L_B\ge 0$, such that 
\begin{eqnarray}
|A(\phi_1)-A(\phi_2)|\le L_B\|\phi_1-\phi_2\|,\;{\rm for\;all\;} \phi_1,\phi_2\in B. 
\nonumber
\end{eqnarray}
\end{definition}
\noindent
We identify $\R^{n\times1}$ and $\R^n$. 
\begin{remark}\rm
Note the $C$-norm on the right hand side. For $\R^n$-valued functionals the property is stronger than the corresponding local property in \cite{Walther2}, which would be implied by (S2-S3), and stronger than the property being {\it almost locally Lipschitz} in \cite{Mallet}. In \cite{Walther} (Lb) is defined for 
$\R^n$-valued functionals and used to show compactness of the maps  $t\_> S(t,\cdot)$ for $t\ge h$.  
\end{remark}

\noindent
If for $f$ there exists a local semiflow let us denote by 
\[
T_\phi:=\{x_t^\phi:\;t\in I_\phi\}\subset C^1
\]
 the trajectory of $\phi\in X$. We denote by $\overline A$ the closure of a set $A$. Our
 main result for general differential equations with state dependent delay can then be
 formulated as
 \begin{theorem}\label{theo1}
Suppose that there exists a local semiflow and that $f$ is (Lb). Let 
$\phi\in X$ be such that $\overline T_\phi\subset U$ and $t_\phi<\infty$. Then for all $C^1$-closed and $C^1$-bounded $L\subset U$ there exists
some $t_L<t_\phi$, such that $x_t\notin L$ for all $t\in[t_L,t_\phi)$. 
\end{theorem}
\noindent
A consequence is a sufficient criterion for global existence:
\begin{corol}\label{cor1}
Suppose that there exists a local semiflow and that $f$ is (Lb). Let $\phi\in X$
be such that $T_\phi$ is $C^1$-bounded and $\overline T_\phi\subset U$,
then $t_\phi=\infty$, i.e., the solution is global. 
\end{corol}
 \noindent
 To prove this, simply assume that $t_\phi<\infty$ and apply Theorem \ref{theo1} to
$L:=\overline T_\phi$, which produces the contradiction that $x_{t_L}\notin L=\overline T_\phi$ for some $t_L<t_\phi$. 
%
\subsection{Variation of constants formula}\label{ss2}
In what follows we study the equation (\ref{eq69}), where $A:U\_> \R^{n\times n}$ is diagonal-matrix-valued, i.e., 
$(A)_{ij}:=\de_{ij}a_i$, 
\begin{eqnarray}
\de_{ij}:=\begin{cases}
1,&\;i=j,
\nonumber\\
0,&\;i\neq j, 
\end{cases}
\nonumber
\end{eqnarray}
for $1\le i,j\le n$ is the Kronecker-delta
and $a,b:U\_> \R^n$ are given. Note that for $\phi\in U$, $x\in\R^n$ and
$\{e_1,...,e_n\}$ the standard basis of $\R^n$, we have
$A(\phi)x=\sum_{i=1}^ne_ia_i(\phi)x_i$.
If we define
\begin{eqnarray}
f(\phi):=A(\phi)\phi(0)+b(\phi) =\sum_{i=1}^ne_ia_i(\phi)\phi_i(0)+b(\phi),
\nonumber
\end{eqnarray}
the DDE $x'(t)=f(x_t)$ becomes (\ref{eq69}). For this type of equations  we 
show a variation of constants formula that  will be useful for showing $C^1$-boundedness and nonnegativity of the trajectory.
\begin{lemma}\label{lem21}
{\rm (Variation of constants formula)} Suppose that for $\phi\in U$ there exists a
solution $x$ on $[-h,t_*)$ of (\ref{eq69}). Then $x$ satisfies
\begin{eqnarray}
x(t)=e^{\int_0^tA(x_s)ds}\left[\phi(0)+\int_0^te^{-\int_0^sA(x_\s)d\s}b(x_s)ds\right],\;\;t\in[0,t_*)
\label{eq70}
\end{eqnarray}
with $(e^{\int_0^tA(x_s)ds})_{ij}=\de_{ij}e^{\int_0^ta_{i}(x_s)ds}$.
\end{lemma}
\noindent
To guarantee that $\overline T_\phi\subset U$ one should make further assumptions on $a
$, $b$ and $U$, so to keep some generality, in the present setting we keep 
$\overline T_\phi\subset U$ as an assumption. 
We now have $X=\{\phi\in U:\;\phi'(0)=A(\phi)\phi(0)+b(\phi)\}$.
\begin{theorem}\label{theo3}
Suppose that $X\neq\emptyset$ and $a$ and $b$ fulfill the smoothness condition (S), then there exists a local semiflow. Let now additionally $a$ and $b$ be (Lb),
map $C^1$-bounded subsets of $U$ on bounded sets and for some $\phi\in X$ let $a$, $b$ and $x$ fulfill the boundedness property
\begin{itemize}
\item[(B)]
there exist $K_1,K_2\ge 0$, such that
\begin{eqnarray}
|a(x_s^\phi)|\le K_1,\;|b(x_s^\phi)|\le e^{K_2s},\;{\rm for \;all} \;s\in(0,t_\phi).
\nonumber
\end{eqnarray}
\end{itemize}
If moreover $\overline T_\phi\subset U$, then 
$t_\phi=\infty$, i.e., the solution is global. 
\end{theorem}
The proof that we will give can be outlined as follows:  The condition that $a$ and $b$ fulfill (S) implies that $f$ fulfills (S) and Theorem \ref{theo2} can be applied. The property
 (Lb) and the condition that $a$ and $b$ map $C^1$-bounded subsets of $U$ on bounded sets imply that $f$ is (Lb). (B) will be used to estimate the trajectory with the variation of constants formula and Corollary \ref{cor1} can be applied.
%
\subsection{A DDE describing stem cell maturation}\label{ss3}
In the following, we show how the assumed conditions for $a$ and $b$ can be satisfied if $a$ and $b$ are specified such that (\ref{eq69}) describes the maturation process of stem cells as modeled in the introduction.
Let $I\subset \mathbb{R}$, $I\neq 0$ be open. We will use the following notation and definitions:
\begin{eqnarray}
&&C[a,b]:=C([a,b],\R),\;\;C^1[a,b]:=C^1([a,b],\R),\;\;M:=C^1([-h,0],I),
\nonumber\\
&&U:=C^1[-h,0]\times M,\;\;\mathbb{R}^m_+:=\{x\in\mathbb{R}^m:\:\;x_i\ge 0,\;i=1,...,m\},
\;m\in\mathbb N.
\nonumber
\end{eqnarray}
Hence, also $M$ and $U$ are open and $\mathbb{R}^m_+$ is the nonnegative cone of $
\mathbb{R}^m$. We now focus on $n=2$ and on delays $\t$ that are allowed to depend on exactly the second of the two components. Suppose that
\begin{eqnarray}
&&\b:I\longrightarrow\mathbb{R}_+,\;\;q:I\longrightarrow\mathbb{R}
,\;\;\t:M\longrightarrow[0,h],\;\;\cF:M\longrightarrow \mathbb{R}_+,\;\;\mu\ge 0.
\nonumber\\
\label{eq75}
\end{eqnarray}
Note that since $\t(M)\subset[0,h]$, for $\vi\in C^1[-h,0]$ and $\psi\in M$ the evaluation 
$\vi(-\t(\psi))$ is well-defined and $\psi(-\t(\psi))\in I$ holds. Then we define
\begin{eqnarray}
a(\vi,\psi)&:=&(q(\psi(0)),-\mu)
\label{eq4}\\
b(\vi,\psi)&:=&(0,\b(\psi(-\t(\psi)))\vi(-\t(\psi))\cF(\psi)).
\label{eq71}
\end{eqnarray}
It follows that 
\begin{eqnarray}
f(\vi,\psi)=
\left(
\begin{array}{c}
q(\psi(0))\vi(0)
\\
\b(\psi(-\t(\psi)))\vi(-\t(\psi))\cF(\psi)-\mu\psi(0)
\end{array}
\right).
\label{eq73}
\end{eqnarray}
If we define $x=(w,v)$ the equation $x'(t)=f(x_t)$ yields (\ref{eq68}--\ref{eq11}).
We first guarantee that $X\neq\emptyset$. We assume that
$0\in I$, as this guarantees that $0\in M$ and $0\in U$. Then, as
$f(0)=0$, we also have that $0\in X$, hence $X\neq\emptyset$. Similarly one can show
that if there exists a nontrivial equilibrium, the corresponding constant function also lies in $X$. We rewrite the variation of constants formula for a given $\phi=(\vi,\psi)\in X$ and for a given solution $(w,v)$ for $t\in[0,t_\phi)$ as
\begin{eqnarray}
w(t)&=&\vi(0)e^{\int_0^tq(v(s))ds},
\label{eq1}
\\
v(t)&=&e^{-\mu t}[\psi(0)+\int_0^t e^{\mu s}\b(v(s-\t(v_s)))\cF(v_s)w(s-\t(v_s))ds]. 
\nonumber\\
\label{eq3}
\end{eqnarray}
To get global existence we would like to guarantee that the closure
of the trajectory lies in $U=C^1[-h,0]\times M$. This can be done if we allow for a large 
range $I$ for the functions in $M$. This leads us to assume that $I=(R_-,\infty)$, for some 
$R_-<0$. 
\begin{theorem}\label{theo4}
Suppose that $I\subset \R$ open, $0\in I$ and $M=C^1([-h,0],I)$. Let $\b$, $q$, $\t$, $\cF$
 and $\mu$ be as in (\ref{eq75}), suppose that
 $\b$ and $q$ are continuously differentiable and $\t$ and $\cF$ satisfy {\rm (S)}. Then
the following properties hold:
\begin{itemize}
\item[(a)] For the DDE describing stem cell maturation (\ref{eq68}--\ref{eq11}) there exists a local semiflow $S$ on $\Ome$ in the sense of Theorem \ref{theo2}.
\item[(b)] If additionally $I=(R_-,\infty)$ for some $R_-<0$, $\t$ and $\cF$ are (Lb), $\cF$ is bounded, $\b$ and $q$ are bounded and Lipschitz on bounded sets, then for  $(t,\phi)\in\Ome$ and $\phi\in C^1([-h,0],\R_+^2)$ one has $S(t,\phi)\in C^1([-h,0],\R_+^2)$ and
$t_\phi=\infty$, i.e., nonnegative initial conditions yield nonnegative global solutions. 
\end{itemize}
\end{theorem}
%
%
\subsection{Specification of $\t$ and $\cF$}\label{ss4}
We denote open balls by $B(x_0,b):=\{x\in\R:\;|x-x_0|<b\}$ for some $x_0\in\mathbb R$ and some $b>0$. 
The following result is a corollary of the Picard-Lindel\"of theorem
and we shall refer to it as the Picard-Lindel\"of theorem. It is our specification of the delay,
 see Figure \ref{f1} for visualization. 
\begin{theorem} {\bf (Picard-Lindel\"of theorem)}\label{theo7}
Let $I\subset \R$ be open, $I\neq\emptyset$.
Suppose that there exist real positive numbers $x_1$, $x_2$,  $b$, $K$ and $\e$ and a function $g$ with the following properties
\begin{itemize}
\item[(i)] $g:\oB\times I\_>\R$ is continuous,  
\item[(ii)] $g$ is uniformly Lipschitz with constant $\frac{K}{b}$ in the first argument,
\item[(iii)]$\e\le g(x,y)\le K$ for $(x,y)\in\oB\times I$ and $x_2-x_1\in(0,\frac{b}{K}\e)$.
\end{itemize}
Then for all $\psi\in C^1([-\frac{b}{K},0],I)$ there exists a unique 
$y(\psi)=y(\cdot,\psi)\in C^1[0,\frac{b}{K}]$
with $y([0,\frac{b}{K}],\psi)\subset\oB$ and
\begin{eqnarray}
\begin{array}{ccc}
y'(s)&=&-g(y(s),\psi(-s)),\;\;s>0,
\\
y(0)&=&x_2.
\end{array}
\label{eq5}
\end{eqnarray}
 Moreover there exists a unique $\t=\t(\psi)\in(0,\frac{b}{K})$, such that $y(\t,\psi)=x_1$.
\end{theorem}
Note that condition (iii) implies that $x_1\in \oB$.
\begin{remark}\rm
Note that as $g$ is continuous, for $a\in(0,h)$ the function
\begin{eqnarray}
\oB\times[0,a]&\longrightarrow&\mathbb{R}
\nonumber\\
(y,s)&\longmapsto&-g(y,\psi(-s))
\nonumber
\end{eqnarray}
is continuous on a compact set and thus bounded. Hence, for a fixed $\psi$ the 
boundedness
of $g$ does not have to be assumed. We assume it nevertheless as we would like to have 
an interval of existence that is uniform for all $\psi$. 
\end{remark}
\noindent
 To achieve that $\t$ fulfills (S), we sharpen the assumptions of the previous theorem. 
More precisely, we will assume that $g$ satisfies property  (G): There exist numbers $x_1, x_2,  b, K, \e\in\R$ and an open interval $J$ with
\begin{itemize}
\item[(G1)] $\oB\subset J$ and ($g$ can be extended such that) $g:J\times I\_>\R$ is $C^1$,  
\item[(G2)] $|\d_1 g (x,y)|<\frac{K}{b}$ for all $(x,y)\in\oB\times I$, 
\item[(G3)] $0<\e\le g(x,y)\le K$ for $(x,y)\in\oB\times I$ and $x_2-x_1\in(0,\frac{b}{K}\e)$.
\end{itemize}
We specify the population growth factor for a given maturation rate $d$, using the 
ingredients $\t$ and $y$ obtained via the Picard-Lindel\"of theorem, as
\[
\cF(\psi)=e^{\int_0^{\t(\psi)}d(y(s,\psi),\psi(-s))ds}.
\]
We now formulate the main result of this paper. 
\begin{theorem}\label{theo5}
Let $I\subset \R$, $I\neq \emptyset$ be open. 
Suppose that $g$ satisfies (G). Let 
$h:=\frac{b}{K}$. Then the following statements hold:
\begin{itemize}
\item[(a)] By the Picard-Lindel\"of theorem for any 
$\psi\in M=C^1([-h,0],I)$ there exists a unique $y=y(\cdot,\psi)\in C^1[0,h]$ with 
$y([0,h],\psi)\subset \oB$ and a unique $\t=\t(\psi)\in[0,h]$, such that $y(\t,\psi)=x_1$. Moreover $\t$ satisfies (S).
\item[(b)] Let additionally $d:J\times I\_>\R$ be $C^1$, then also $\cF$ satisfies (S).
\item[(c)] Suppose that moreover $0\in I$,  $\b$, $q$ and $\mu$ are as in (\ref{eq75})
and $f$ is as in (\ref{eq73}). 
Then $f$ induces a local semiflow $S$ on $
\Omega$. 
\item[(d)] Suppose that additionally $I=(R_-,\infty)$, and that  the sets
\[
d(\oB\times I),\;\d_2 g (\oB\times A),\;{\rm and}\;\d_id(\oB\times A),\;i=1,2,
\]
are bounded, whenever $A\subset I$ is bounded and that $\b$ and $q$ are bounded and Lipschitz on bounded sets.
Let $(t,\phi)\in\Ome$ and $\phi\in C^1([-h,0],\R_+^2)$, then 
$S(t,\phi)\in C^1([-h,0],\R_+^2)$ and $t_\phi=\infty$, i.e., nonnegative initial conditions yield
nonnegative global solutions. 
\end{itemize}
\end{theorem}
%
%
\section{Proofs of Section \ref{s1}}\label{s2}
\subsection{Differential equations with state dependent delay - proofs of \ref{ss1}}\label{S1}
%
To prove that in case of a finite existence interval closed and 
bounded sets are left for good (Theorem \ref{theo1}) we will apply the existing result
that compact sets are left for good (Lemma \ref{lem3}). A further useful tool is the following
result on $C^1$-compactness that is a straightforward corollary (which we state without proof) of the Arzela-Ascoli theorem. 
For $A\subset C^1$, we denote by 
$A':=\{f':f\in A\}\subset C$ the set of derivatives of $A$. 
\begin{lemma}\label{lem18}
If $A\subset C^1$ is $C^1$-bounded and $A$ and $A'$ are equicontinuous, then $A$ is
relatively $C^1$-compact. 
\end{lemma}

 \bigskip
\noindent
{\bf Proof of Theorem \ref{theo1}.}  Let $\phi\in X$ and $L\subset U$ be $C^1$-closed and $C^1$-bounded. Choose $r$, such that $\psi\notin L$, 
whenever $\|\psi\|\ge r$ or $\|\psi'\|\ge r$. We first consider the case that $x$ is 
unbounded. By its continuity $x$ is bounded
on $[-h,t_\phi-h]$. From the unboundedness of $x$ it then follows that there exists some
$t_N\in[t_\phi-h, t_\phi)$, such that $|x(t_N)|\ge r$. Let $t\in [t_N,t_\phi)$. 
Then
\[
\|x_t\|=\sup_{\th\in[-h,0]}|x(t+\th)|\ge |x(t_N)|\ge r. 
\]
Thus $x_t\notin L$. Hence, the conclusion of the theorem holds. Now we assume that
$x$ is bounded on $[-h,t_\phi)$. Choose $M_1>0$, such that $\|x_t\|\le M_1$ for all $t\in I_\phi$. Thus, $T_\phi$ is $C$-bounded. Next, we consider the case that $x'$ is unbounded. By its continuity $x'$
is bounded on $[-h,t_\phi-h]$. From the unboundedness of $x'$ it then follows that there 
exists some $t_N\in[t_\phi-h, t_\phi)$, such that $|x'(t_N)|\ge r$. Let $t\in [t_N,t_\phi)$. 
Then 
\[
\|x'_t\|=\sup_{\th\in[-h,0]}|x'(t+\th)|\ge |x'(t_N)|\ge r. 
\]
Thus $x_t\notin L$ and again the conclusion of the theorem holds. 
Hence, we should consider the case that there is some $M_2>0$, such that 
$\|x_t'\|\le M_2$ for all $t\in  I_\phi$. It follows that $T_\phi'$ is $C$-bounded and
thus $T_\phi$ is $C^1$-bounded. From the boundedness of $x'$ we can conclude
that there exists some $M_3>0$, such that
\[
|x'(t)|\le M_3\;\;{\rm for\;all\;} t\in I_\phi.
\]
Hence for $t_1\ge t_2\ge 0$ 
 one has
\[
|x(t_1)-x(t_2)|\le\int_{t_2}^{t_1}|x'(t)|dt\le M_3|t_1-t_2|
\]
and thus $x$ is uniformly continuous on $I_\phi$ and thus also on $[-h,t_\phi)$. It 
follows that $T_\phi$ is equicontinuous. 
As $f$ is (Lb) and $x'$ is bounded, we know that
for $t\ge s>0$
\begin{eqnarray}
|x'(t)-x'(s)|&=&|f(x_t)-f(x_s)|\le L_{T_\phi}\|x_t-x_s\|
\nonumber\\
&=&L_{T_\phi}\sup_{\th\in[-h,0]}|x(t+\th)-x(s+\th)|
\nonumber\\
&=& L_{T_\phi}\sup_{\th\in[-h,0]}\int_{s+\th}^{t+\th}|x'(\s)|d\s\le L_{T_\phi} b|t-s|
\nonumber
\end{eqnarray}
for some constant $b$. Thus $x'$ is uniformly continuous on $I_\phi$ and thus also on
$\overline I_\phi$. Hence, $T_\phi'$ is equicontinuous. We have shown that $T_\phi$
satisfies the assumptions of Lemma \ref{lem18}. 
Hence,  $\overline T_\phi$ is $C^1$-compact. The assumption $\overline T_\phi\subset U$ 
and the continuity of $f$ imply that $\overline T_\phi\subset X$. Then
Lemma \ref{lem3} implies the existence of some $t_{\overline T_\phi}$, such that 
$x_t\notin\overline T_\phi$ for all $t\in[t_{\overline T_\phi}, t_\phi)$, which is a 
contradiction. 
\qed\bigskip

\subsection{Variation of constants formula - proofs of \ref{ss2}}
\label{ss5}
To show existence of a local semiflow it remains to guarantee that 
$f(\phi)=A(\phi)\phi(0)+b(\phi)$ satisfies (S):
\begin{lemma}
Suppose that $a$ and $b$ satisfy (S), then so does $f$ and 
\begin{eqnarray}
Df(\phi)\chi&=&A(\phi)\chi(0)+\phi(0)DA(\phi)\chi+Db(\phi)\chi
\nonumber\\
(DA(\phi)\chi)_{ij}&=&\de_{ij}Da_i(\phi)\chi_i
\nonumber
\end{eqnarray}
for $\phi\in U$, $\chi\in C^1$. The extension $D_ef$ required in (S) is obtained by replacing $DA(\phi)$
and $Db(\phi)$ by $D_eA$ and $D_eb$, where 
$(D_eA(\phi)\chi)_{ij}=\de_{ij}D_ea_i(\phi)\chi_i$. 
\end{lemma}
\Proof Define $z:U\_>\R^n; \phi\longmapsto\phi(0)$. Then $z$ is $C^1$ with
$Dz(\phi)\chi=\chi(0)$ and this formula can also define an extension $D_ez(\phi)$ to
$C$. One easily checks that 
\[
U\times C\_> \R^n;\;(\phi,z)\_> D_ez(\phi)\chi
\]
is continuous. Thus, $z$ satisfies (S). As so do $a$ and $b$ it follows
from the sum- and product rules of differentiation and continuity that $f$ satisfies (S) and has the derivatives and extensions as claimed. 
\qed\bigskip

To show global existence we first guarantee that $f$ is (Lb). 
The proof of the following result is straightforward and we omit it. 
\begin{lemma}\label{lem19}
If $f_1,f_2:U\_>\R^{n\times m}$ are (Lb), then so is $f_1+f_2$.
If additionally $f_1$ and $f_2$ map $C^1$-bounded subsets of $U$ on bounded
sets, then $f_1\cdot f_2$ is (Lb).
\end{lemma}
A consequence is
\begin{lemma}
If $a$ and $b$ are (Lb) and map $C^1$-bounded subsets of $U$ on bounded
sets then $f$ is (Lb). 
\end{lemma}

\Proof By definition of $A$ and the hypothesis it follows that $A$ is (Lb)
and  maps $C^1$-bounded subsets of $U$  on bounded sets.  
Then, by Lemma \ref{lem19}, $f$ is (Lb). 
\qed\bigskip

Next, to estimate the trajectory we show that the variation of constants formula holds. 

\bigskip
\noindent
{\bf Proof of Lemma \ref{lem21}.} Define $y(t):=e^{-\int_0^tA(x_s)ds}x(t)$, $t\in[0,t_*)$. By the solution properties of
$x$ and the product rule of differentiation we get that $y$ is $C^1$ on $[0,t_*)$. For
$t\in[0,t_*)$ one has
\begin{eqnarray}
y'(t)&=&-A(x_t)e^{-\int_0^tA(x_s)ds}x(t)+e^{-\int_0^tA(x_s)ds}x'(t)
\nonumber\\
&=&-A(x_t)e^{-\int_0^tA(x_s)ds}x(t)+e^{-\int_0^tA(x_s)ds}[A(x_t)x(t)+b(x_t)]
\nonumber\\
&=&e^{-\int_0^tA(x_s)ds}b(x_t).
\nonumber
\end{eqnarray}
It follows that 
\begin{eqnarray}
y(t)=\int_0^te^{-\int_0^sA(x_\s)d\s}b(x_s)ds+\phi(0),\;\;t\in[0,t_*).
\nonumber
\end{eqnarray}
We use the definition of $y$ and get the desired result. 
\qed\bigskip

\begin{lemma}\label{lem20}
Suppose that there exists a local semiflow, that for some $\phi\in X$ one has $t_\phi<\infty$ and that $a$, $b$ and $x$ fulfill the boundedness property (B), 
 then $T_\phi$ is $C^1$-bounded. 
\end{lemma}

\Proof With the variation of constants formula we can estimate
\begin{eqnarray}
|x(t)|&\le& \|\phi\|^1|e^{\int_0^tA(x_s)ds}|+
|e^{\int_0^tA(x_s)ds}\int_0^te^{-\int_0^sA(x_\s)d\s}b(x_s)ds|
\nonumber\\
&\le& \|\phi\|^1e^{\int_0^t|a(x_s)|ds}+
\int_0^te^{\int_s^t|a(x_\s)|d\s}|b(x_s)|ds.
\nonumber
\end{eqnarray}
With the assumed boundedness properties of $a$ and $b$ and since $t_\phi<\infty$
it follows that $\{x(t):t\in I_\phi\}$ is bounded. Next, we can estimate the expression
\[
|x'(t)|=|A(x_t)x(t)+b(x_t)|
\]
using the conditions imposed on $a$ and $b$ and the boundedness of $x$. As 
$\phi\in U\subset C^1$, we also have boundedness of $x$ and $x'$ on $[-h,0]$ and it
follows that $T_\phi$ is $C^1$-bounded. 
\qed\bigskip

Now we are ready to apply the results for general differential equations with state 
dependent delay to show local and global existence for (\ref{eq69}). 

\bigskip
\noindent{\bf Proof of Theorem \ref{theo3}.} As $a$ and $b$ satisfy (S), so does $f$
and the first statement follows by Theorem \ref{theo2}. Next, note that under the conditions
of the theorem we have guaranteed that $f$ is (Lb). Suppose now that $\phi\in X$ with $t_\phi<\infty$. Then, by Lemma \ref{lem20}, $T_\phi$ is $C^1$-bounded. As 
$\overline T_\phi\subset U$ is supposed, we get $t_\phi=\infty$ with Corollary \ref{cor1}. \qed\bigskip

\subsection{A DDE describing stem cell maturation - proofs of \ref{ss3}}\label{S2}
\subsubsection{Noncontinuable solutions}
We guarantee that the conditions of Theorem \ref{theo3} hold. 
%
%
First we guarantee that $a$ and $b$ as defined in (\ref{eq4}--\ref{eq71}) satisfy (S). 
\begin{prop}\label{prop3} Let $\b$, $q$, $\t$, $\cF$ and $\mu$ be as
in (\ref{eq75}), suppose that $\b$ and $q$ are continuously differentiable and that $\t$
and $\cF$ satisfy (S) on $M$. Then $a$ and $b$ satisfy
(S) on $U$. For $(\vi,\psi)\in U$ and $\chi,\xi\in C^1[-h,0]$, one has
\begin{eqnarray}
&&Da(\vi,\psi)(\chi,\xi)=(q'(\psi(0))\xi(0),0),
\nonumber\\
&&Db(\vi,\psi)(\chi,\xi)=D(b_1,b_2)(\vi,\psi)(\chi,\xi),\;\;{\rm where}
\nonumber\\
&&Db_1(\vi,\psi)(\chi,\xi)=0,
\nonumber\\
&&Db_2(\vi,\psi)(\chi,\xi)=
\nonumber\\
&&\vi(-\t(\psi))\cF(\psi)\b'(\psi(-\t(\psi)))[\xi(-\t(\psi))-
\psi'(-\t(\psi))D\t(\psi)\xi]
\nonumber\\
&&+\cF(\psi)\b(\psi(-\t(\psi)))[\chi(-\t(\psi))-\vi'(-\t(\psi))D\t(\psi)\xi]
\nonumber\\
&&+\b(\psi(-\t(\psi)))\vi(-\t(\psi))D\cF(\psi)\xi.
\nonumber
\end{eqnarray}
One gets the extensions $D_ea$ and $D_eb$ if one replaces $D\t(\psi)$ and $D\cF(\psi)$ by the
respective extensions. 
\end{prop}
To prove these results we introduce the evaluation operator, 
\begin{eqnarray}
Ev:M\times[-h,0]&\longrightarrow&\mathbb{R}
\nonumber\\
(\phi,s)&\longmapsto&\phi(s).
\nonumber
\end{eqnarray}
Here $M$ can be an arbitrary open subset of $C^1[-h,0]$. It is shown in \cite[p 481]{Walther} that $Ev$ is continuously differentiable with
\begin{eqnarray}
D_1Ev(\phi,s)\chi=Ev(\chi,s),\;\;D_2Ev(\phi,s)1=\phi'(s).
\nonumber
\end{eqnarray}
Then one can use the identity
\begin{eqnarray}
Ev\circ(id\times-\t)(\vi,\psi)=\vi(-\t(\psi))
\nonumber
\end{eqnarray}
and show continuous differentiability of $a$ and $b$ with the chain rule. The rest of the proof is standard. 
The result was shown in \cite{Alarcon} for the case where $I$ is a bounded open
interval. In summary we have 
\begin{prop}\label{prop2}
Suppose that $\b$ and $q$ are continuously differentiable, $\t$ and $\cF$
satisfy (S) and $X\neq 0$. 
Then for $a$ and $b$ defined by (\ref{eq4}--\ref{eq71}) there exists a local semiflow. 
\end{prop}
\Proof This is a combination of Proposition \ref{prop3} and Theorem \ref{theo3}. 
\qed\bigskip

\subsubsection{Nonnegativity}
With the variation of constants formula it is easy to see that the semiflow maps nonnegative times and initial functions to 
nonnegative functions. 
\begin{lemma}\label{lem23}
Suppose that there exists a local semiflow $S$ on $\Ome$ as defined in Theorem
\ref{theo2},
 let $(t,\phi)\in\Ome$ and suppose that $\phi\in C^1([-h,0],\mathbb{R}_+^2)$, then 
$S(t,\phi)\in  C^1([-h,0],\mathbb{R}_+^2)$. 
\end{lemma}
\Proof
Let $\phi=(\vi,\psi)\in C^1([-h,0],\R^2_+)$, $(t,\phi)\in\Ome$ and set $x^\phi=(w,v)$. Then 
$(w,v)$ satisfies (\ref{eq1}--\ref{eq3}) by Lemma \ref{lem21}. By (\ref{eq1}) one has $w(t)\ge 0$. Nonnegativity of $\psi$, $\b$, $\cF$ and $w$ then imply via
(\ref{eq3}) that $v(t)\ge 0$. 
Hence, $x^\phi\in C^1([-h,T],\mathbb{R}_+^2)$ for all $T\in I_\phi$. Thus
$S(t,\phi)=x_t^\phi\in C^1([-h,0],\mathbb{R}_+^2)$. 
\qed\bigskip

\subsubsection{Global solutions}
To get global existence via Theorem \ref{theo3}, we guarantee that the assumptions of the 
theorem hold. The boundedness and Lipschitz properties of $a$ and $b$ can be shown
in a standard way if we suppose that the new model ingredients satisfy 
corresponding properties:
\begin{lemma}\label{lem22}
Suppose that $\t$ and $\cF$ are (Lb), $\cF$ is bounded, $\b$ and $q$ are
bounded on bounded sets and Lipschitz on bounded sets, then $a$ and $b$ defined in (\ref{eq4}--\ref{eq71}) are (Lb) and bounded on $C^1$-bounded sets. 
\end{lemma}
\Proof Let $B\subset U$ be $C^1$-bounded. First note that $\{\psi(0):\;(\vi,\psi)\in B\}$
is a bounded subset of $\R$.  As $q$ is Lipschitz on bounded sets, it follows that $a$ is (Lb).  
As $q$ is bounded on bounded sets it also follows that $a$ is bounded on $C^1$-bounded sets. To show
the statements for $b$ we should show them for $b_2$, the second component of $b$. 
First note that 
$\{\vi(-\t(\psi)):(\vi,\psi)\in B\}$ and $\{\psi(-\t(\psi)):(\vi,\psi)\in B\}$ are bounded subsets of
$\R$. Next, one can show with the mean value theorem, using differentiability of $\vi$
and $\psi$ and $C^1$-boundedness of $B$, that the maps from 
$U$ to $\R$ given by
\[
(\vi,\psi)\longmapsto\vi(-\t(\psi)){\rm\;and\;}(\vi,\psi)\longmapsto\psi(-\t(\psi))
\]
are (Lb) and bounded on $B$. Then $(\vi,\psi)\longmapsto\vi(-\t(\psi))\cF(\psi)$
is (Lb) by Lemma \ref{lem19} and maps 
$C^1$-bounded sets on bounded sets. Hence $b_2$ also has these properties. 
\qed\bigskip

\begin{lemma}\label{lem4}
Suppose that $\b$, $\cF$ and $q$ are bounded and that there exists a 
local semiflow. Then $a$ and $b$ satisfy the boundedness property (B). 
\end{lemma}

\Proof 
Let $\phi\in X$ and denote by $x=x^\phi=(w,v)$ the noncontinuable solution.  
The boundedness of $q$ implies the boundedness property of $a$ that is required in (B). It remains to show that $b_2$ fulfills the corresponding property. For $t\in(0,t_\phi)$ one has
\begin{eqnarray}
|b_2(x)|=|\b(v(t-\t(v_t)))w(t-\t(v_t))\cF(v_t)|\le K_1|w(t-\t(v_t))|
\label{eq72}
\nonumber
\end{eqnarray}
for some $K_1\ge 0$. From the variation of constants formula and the boundedness of
$q$ we derive an exponential bound for $w$. If we continue estimating $|b_2(x)|$
with this exponential bound, the exponential boundedness property for $b_2$ follows. 
\qed\bigskip

\bigskip
\noindent{\bf Proof of Theorem \ref{theo4}.} (a) follows by Proposition \ref{prop2} and the
nonnegativity statement in (b) follows by Lemma \ref{lem23}.
Global existence follows by Theorem \ref{theo3} and the previous results if we show that 
$\overline T_\phi\subset U$. As $C^1([-h,0],\R_+^2)$ is $C^1$-closed and as, by our 
nonnegativity results and the assumption $I=(R_-,\infty)$, we have 
$T_\phi\subset C^1([-h,0],\R_+^2)\subset U$, we have indeed $\overline T_\phi\subset U$.
\qed\bigskip

%
%
\subsection{Specification of $\t$ and $\cF$ - proofs of \ref{ss4}}
In this section we guarantee that $\t$ and $\cF$ satisfy the required smoothness 
properties. 
%
\subsubsection{Noncontinuable solutions}
%
%
We first show that the ODE induced by $g$ has a solution $y$ and then 
use this solution to define the delay. 

\bigskip
\noindent{\bf Proof of Theorem \ref{theo7} (Picard-Lindel\"of theorem).}
Choose $\psi\in C^1([-\frac{b}{K},0],I)$ and define 
\[
f:[0,\frac{b}{K}]\times\oB\_>\R; \;f(s,y):=-g(y,\psi(-s)).
\] Then $K$ is a bound for
$|f(s,y)|$ and $f$ is continuous. Moreover, $f$ is uniformly
Lipschitz with constant $K/b$ in the second argument. Then, the Picard-Lindel\"of
theorem (see e.g. Theorem II.1.1 in \cite{Hartman}) implies the first statement. 
Next, integration of the ODE yields
\begin{eqnarray}
y(t,\psi)=x_2-\int_0^tg(y(s,\psi),\psi(-s))ds.
\label{eq10}
\end{eqnarray}
Then the second statement follows from the intermediate value theorem applied to 
$t\longmapsto y(t,\psi)$ when using condition (iii). 
\qed\bigskip

\noindent
For the following exposition we assume that $g$ satisfies (G). If we use that $y(\t,\psi)=x_1$, we get that 
\begin{eqnarray}
x_2-x_1=\int_0^{\t(\psi)}g(y(s,\psi),\psi(-s))ds
\label{eq9}
\end{eqnarray}
and, as we will see, $D\t(\psi)$ can be defined by differentiating this equation
with respect to $\psi$. This involves differentiating $y(\cdot,\psi)$ with respect to $\psi$. 
To apply the implicit function theorem, we define a map on an open domain via
\begin{eqnarray}
G: M\times C^1([0,h],J)&\longrightarrow& C^1[0,h]
\nonumber\\
(\psi,y)&\longmapsto&G(\psi,y)
\nonumber\\
&&G(\psi,y)(t):=y(t)+\int_0^tg(y(s),\psi(-s))ds.
\nonumber
\end{eqnarray}
By the Picard-Lindel\"of theorem for any $\psi\in M$ the equation 
\begin{eqnarray}
G(\psi,y)=x_2
\label{eq44}
\end{eqnarray}
as an identity in $C^1[0,h]$ has a solution 
\[
y(\psi)=y(\cdot,\psi)\in C^1([0,h],\oB)\subset C^1([0,h],J).
\] 
To show that $G$ is $C^1$ we use the following
technical result.
\begin{lemma}\label{lem1} Let $E$ be an open interval, 
$a,b\in\R$, $0\le c\le 1$. For any $\ops\in M_E:=C^1([a,b],E)$ there exists some
$d=d(\ops)$, such that for all $\de\in(0,d]$
\begin{eqnarray}
&&A_{\de}
:=
\nonumber\\
&&\{\ops(s)+\th(\psi(s)-\ops(s)):s\in[a,b],\th\in[c,1],\psi\in M_E,\|\psi-\ops\|^1\le\de\}
\nonumber
\end{eqnarray}
is a compact subset of $E$. 
\end{lemma}
\Proof Let $\ops\in M_E$. By the openness of $E$ and the compactness of $[a,b]$ with
\begin{eqnarray}
d(\ops):=\begin{cases}
\frac{1}{2} {\rm dist}(\ops([a,b]),\R\backslash E),&\R\backslash E\neq\emptyset,
\\
\frac{1}{2},&\R\backslash E=\emptyset,
\end{cases}
\end{eqnarray}
where ${\rm dist}$ denotes the distance function, we get $d(\ops)\in(0,\infty)$
and $A_\de\subset E$ for all $\de\in(0,d(\ops)]$. Now fix some $\de\in(0,d(\ops)]$.  Then $A_\de$ is bounded. It remains to show that $A_\de$ is closed. Let $(z_n)\in A_\de^\N$, $z\in\R$, 
$z_n\_>z$, $(\psi_n)\in M_E^\N$ with $\|\psi_n-\ops\|^1\le\de$, $(s_n)\in[a,b]^\N$, 
$(\th_n)\in[c,1]^\N$ such that $z_n=\ops(s_n)+\th_n(\psi_n(s_n)-\ops(s_n))$.
First note that $z_n\in[\inf_{s\in[a,b]}\ops(s)-\de,\sup_{s\in[a,b]}\ops(s)+\de]$. 
Then also $z\in[\inf_{s\in[a,b]}\ops(s)-\de,\sup_{s\in[a,b]}\ops(s)+\de]$. 
Next, we can use 
\begin{eqnarray}
&&z\in[\inf_{s\in[a,b]}\ops(s)-\de,\sup_{s\in[a,b]}\ops(s)+\de]
\nonumber\\
&=&[\inf_{s\in[a,b]}\ops(s)-\de,\inf_{s\in[a,b]}\ops(s)]\cup[\inf_{s\in[a,b]}\ops(s),\sup_{s\in[a,b]}\ops(s)]\nonumber\\
&&\cup[\sup_{s\in[a,b]}\ops(s),\sup_{s\in[a,b]}\ops(s)+\de]
\nonumber
\end{eqnarray}
and the intermediate value theorem to show that there exists some $\os\in[a,b]$ with $|\ops(\os)-z|\le \de$.
For such an $\os$ we define $\psi(s):=\ops(s)-(\ops(\os)-z)$ for all $s\in[a,b]$. Then $\psi$ is $C^1$.
Moreover $\|\psi-\ops\|^1=|\ops(\os)-z|\le \de$. It follows that $\psi\in M_E$. 
Then we choose $
\th:=1$
and get $\ops(s)+\th(\psi(s)-\ops(s))=\psi(s)$. Hence $\psi(s)\in A_\de$ for any $s$.
Since $\psi(\os)=z$, also $z\in A_\de$. Thus, $A_\de$ is closed. 
\qed\bigskip

\begin{remark}\label{rem1}\rm
With a similar proof one can show a variant of the previous result for continuous functions,
i.e., the statement with $C^1([a,b],J)$ replaced by $C([a,b],J)$ and $\|\cdot\|^1$ replaced
by $\|\cdot\|$. 
\end{remark}
We can use the previous result to show some convergence properties that will in turn 
help to show differentiability and continuity  properties of several operators. In the following
we denote by $k$ an arbitrary  function for which we show properties that will be applied
to functions of the model and their derivatives. 
\begin{lemma}\label{lem11}
Suppose that $k:J\times I\_>\R$ is $C^1$. Then the following statements hold:
\begin{itemize}
\item[(a)] Fix $\psi\in M$ and $\oz\in C^1([0,h],J)$. Then for $z\in C^1([0,h],J)$
\begin{eqnarray}
\begin{array}{l}
\sup_{s\in[0,h]}|k(z(s),\psi(-s))-k(\oz(s),\psi(-s))
\end{array}
\nonumber \\
\begin{array}{r}
-\d_1k(\oz(s),\psi(-s))(z(s)-\oz(s))|
\end{array}
\nonumber
\end{eqnarray}
is $o(\|z-\oz\|^1)$ as $\|z-\oz\|^1\rightarrow 0$.
\item[(b)] Fix $z\in C^1([0,h],J)$ and $\ops\in M$. Then for $\psi\in M$
\begin{eqnarray}
\begin{array}{l}
\sup_{s\in[0,h]}|k(z(s),\psi(-s))-k(z(s),\ops(-s))
\end{array}
\nonumber \\
\begin{array}{r}
-\d_2k(z(s),\ops(-s))(\psi(s)-\ops(s))|
\end{array}
\nonumber
\end{eqnarray}
is $o(\|\psi-\ops\|^1)$ as $\|\psi-\ops\|^1\rightarrow 0$. 
\end{itemize}
\end{lemma}
\Proof  By the mean value theorem for $s\in[0,h]$ there exists some $\th_s\in[0,1]$, such that 
\begin{eqnarray}
&&|k(z(s),\psi(-s))-k(\oz(s),\psi(-s))-\d_1k(\oz(s),\psi(-s))(z(s)-\oz(s))|
\nonumber\\
&=&|\d_1k(\oz(s)+\th_s(z(s)-\oz(s)),\psi(-s))-\d_1k (\oz(s),\psi(-s))||z(s)-\oz(s)|
\nonumber\\
&\le&\|z-\oz\|^1|\d_1k (\oz(s)+\th_s(z(s)-\oz(s)),\psi(-s))-\d_1k (\oz(s),\psi(-s))|.
\label{eq47}
\end{eqnarray}
In Lemma \ref{lem1} we choose $[a,b]:=[0,h]$, $E:=J$ and $c=0$. Then 
$M_E=C^1([0,h],J)$ and we know that there exists some $d=d(\oz)$ such that
\begin{eqnarray}
A_{\de}
&:=&\{\oz(s)+\th(z(s)-\oz(s)):s\in[0,h],\th\in[0,1],z\in M_E,\|z-\oz\|^1\le\de\}
\nonumber
\end{eqnarray}
is compact for all $\de\in(0,d(\oz)]$. Also $\psi[-h,0]$ is compact. So $\d_1k $ is  uniformly 
continuous on the compact set $A_\de\times\psi[-h,0]$. 
Then (a) follows from (\ref{eq47}). The proof of (b) is similar, so we omit it. 
\qed\bigskip

\begin{lemma}\label{lem12}
Let $k:J\times I\_>\R$ be continuous. Fix $(\oz,\ops)\in C^1([0,h],J)\times M$, then
for $(z,\psi)\in C^1([0,h],J)\times M$ there exists some $d>0$, such that for all 
$\de\in(0,d]$ the set 
\begin{eqnarray}
&&\{(z(s),\psi(-s)):
\nonumber\\
&&s\in[0,h],(z,\psi)\in C^1([0,h],J)\times M,\|(z,\psi)-(\oz,\ops)\|^1\le\de\}
\nonumber
\end{eqnarray}
is compact. Moreover
\begin{eqnarray}
\sup_{s\in[0,h]}|k(z(s),\psi(-s))-k(\oz(s),\ops(-s))|
\nonumber
\end{eqnarray}
tends to zero if $\|z-\oz\|^1$ and $\|\psi-\ops\|^1$ tend to zero.
\end{lemma}

\Proof In Lemma \ref{lem1} choose $[a,b]:=[0,h]$, $E:=J$, $c:=1$, then 
$M_E=C^1([0,h],J)$ and we know that 
\begin{eqnarray}
A_\de=\{z(s):s\in[0,h],z\in C^1([0,h],J),\|z-\oz\|^1\le\de\}
\nonumber
\end{eqnarray}
is compact for all $\de$ in a neighborhood of zero. Similarly,
\begin{eqnarray}
B_\de=\{\psi(-s):s\in[0,h],\psi\in M,\|\psi-\ops\|^1\le\de\}
\nonumber
\end{eqnarray}
is compact for all $\de$ in a neighborhood of zero. Hence, there exists some $d>0$, 
such that $A_\de\times B_\de$ is compact for all $\de\in(0,d]$ and the first statement is 
shown. Thus the stated limit behavior follows from uniform continuity of $k$. 
\qed\bigskip

\begin{remark}\rm
The previous result can be formulated for continuous functions with a similar proof, if
one uses a variant of Lemma \ref{lem1}, see Remark \ref{rem1}. 
\end{remark}
In the following, to save brackets, we write $Ax(t)$ instead of $(Ax)(t)$ to denote 
evaluated functions in the image of a linear operator. By $\|\cdot\|_{op}$ we denote
operator norms. We are now ready to prove
\begin{lemma}\label{lem6}
The operator $G$ is $C^1$ with derivative $DG=(\d_1G ,\d_2 G )$, where for 
$(\psi,y)\in M\times C^1([0,h],J)$
\begin{eqnarray}
&&\d_1G (\psi,y):C^1[-h,0]\_>C^1[0,h],
\nonumber\\
&&
\d_1G (\psi,y)\chi(t)=\int_0^t\d_2 g (y(s),\psi(-s))\chi(-s)ds,
\nonumber\\
&&\d_2 G (\psi,y):C^1[0,h]\_>C^1[0,h],
\nonumber\\
&&\d_2 G (\psi,y)z(t)=z(t)+\int_0^t\d_1 g (y(s),\psi(-s))z(s)ds.
\nonumber
\end{eqnarray}
\end{lemma}

\Proof
It is sufficient to show that the partial derivatives 
$\d_1G (\psi,y)$ and $\d_2 G (\psi,y)$ exist and the maps $(\psi,y)\longmapsto \d_iG(\psi,y)$, $i=1,2$
are continuous. 
First, note that $\d_1G (\psi,y)$ is a well-defined bounded linear operator. 
Then 
\begin{eqnarray}
&&|[G(\psi,y)-G(\ops,y)-\d_1G (\ops,y)(\psi-\ops)](t)|\le
\nonumber\\
&&\int_0^h|g(y(s),\psi(-s))-g(y(s),\ops(-s))-\d_2 g (y(s),\ops(-s))(\psi-\ops)(-s))|ds.
\nonumber
\end{eqnarray}
By Lemma \ref{lem11} (b) the expressions are $o(\|\psi-\ops\|^1)$ as $\|\psi-\ops\|^1\rightarrow 0$ uniformly in $t$
and thus also
\begin{eqnarray}
&&\|G(\psi,y)-G(\ops,y)-\d_1G (\ops,y)(\psi-\ops)\|
\nonumber
\end{eqnarray}
is $o(\|\psi-\ops\|^1)$ as $\|\psi-\ops\|^1\rightarrow 0$. Next, 
\begin{eqnarray}
&&|[G(\psi,y)-G(\ops,y)-\d_1G (\ops,y)(\psi-\ops)]'(s)|
\nonumber\\
&=&|g(y(s),\psi(-s))-g(y(s),\ops(-s))-\d_2 g (y(s),\ops(-s))(\psi(-s)-\ops(-s))|.
\nonumber
\end{eqnarray}
Again, the expression is $o(\|\psi-\ops\|^1)$ as $\|\psi-\ops\|^1\rightarrow 0$ uniformly in $s$.
Hence
\begin{eqnarray}
&&\|[G(\psi,y)-G(\ops,y)-\d_1G (\ops,y)(\psi-\ops)]'\|
\nonumber
\end{eqnarray}
is $o(\|\psi-\ops\|)^1$ as $\|\psi-\ops\|^1\rightarrow 0$. Thus, so is
\begin{eqnarray}
&&\|G(\psi,y)-G(\ops,y)-\d_1G (\ops,y)(\psi-\ops)\|^1.
\nonumber
\end{eqnarray}
This implies partial differentiability in the first argument. Next, 
we show continuity of $(\psi,y)\longmapsto \d_1G (\psi,y)$ on $M\times C^1([0,h],J)$. 
First,
\begin{eqnarray}
&&|[\d_1G (\psi,y)\chi-\d_1G (\ops,\oy)\chi](t)|
\nonumber\\
&\le&
\|\chi\|^1\int_0^h|\d_2 g (y(s),\psi(-s))-\d_2 g (\oy(s),\ops(-s))|ds.
\nonumber
\end{eqnarray}
The integral tends to zero by Lemma \ref{lem12} if $\|y-\oy\|^1$ and $\|\psi-\ops\|^1$
tend to zero. Hence for any $\e>0$, there exists a $\de=\de(\e,\ops,\oy)$, such that 
\begin{eqnarray}
&&\|[\d_1G (\psi,y)-\d_1G (\ops,\oy)]\chi\|\le\frac{\e}{2}\|\chi\|^1
\nonumber
\end{eqnarray}
for $\|\psi-\ops\|^1\le \de$, $\|y-\oy\|\le\de$. Similarly, one can show that there exists a
$\de$, such that 
\begin{eqnarray}
&&\|[\d_1G (\psi,y)\chi-\d_1G (\ops,\oy)\chi]'\|\le\frac{\e}{2}\|\chi\|^1
\nonumber
\end{eqnarray}
for $\|\psi-\ops\|^1\le \de$, $\|y-\oy\|\le\de$ and thus also a $\de$, such that
\begin{eqnarray}
&&\|[\d_1G (\psi,y)-\d_1G (\ops,\oy)]\chi\|^1\le\e\|\chi\|^1
\nonumber
\end{eqnarray}
for $\|\psi-\ops\|^1\le \de$, $\|y-\oy\|\le\de$.
It follows that in the operator norm we get
\begin{eqnarray}
\|\d_1G (\psi,y)-\d_1G (\ops,\oy)\|_{op}\le\e
\nonumber
\end{eqnarray}
for $\|\psi-\ops\|^1\le \de$, $\|y-\oy\|\le\de$ and this implies the desired continuity result. Next, note that also $\d_2 G (\psi,y)$ is a well-defined 
bounded linear operator. 
The proof of partial differentiability and continuity of $(\psi,y)\longmapsto \d_2 G (\psi,y)$ is
similar as in the first argument. Thus $G$ is $C^1$. 
\qed\bigskip

The next result follows from the theory of linear ODE, in particular the variation of
constants formula for ODE. 
\begin{lemma}\label{lem7}
The operator $\d_2 G (\psi,y)$ has the bounded linear inverse
\begin{eqnarray}
\d_2 G (\psi,y)^{-1}&:& C^1[0,h]\longrightarrow C^1[0,h]
\nonumber\\
\d_2 G (\psi,y)^{-1}x(t)
&=&x(0)e^{-\int_0^t\d_1 g (y(s),\psi(-s))ds}
\nonumber\\
&&+\int_0^te^{-\int_s^t\d_1 g (y(\s),\psi(-\s))d\s}x'(s)ds.
\nonumber\\
\label{eq18}
\end{eqnarray}
\end{lemma}
\Proof First, $\d_2 G (\psi,y)z=x$ for $z,x\in C^1[0,h]$
can be written as 
\begin{eqnarray}
z(t)+\int_0^tK(s)z(s)ds=x(t),\;\;K(s):= \d_1 g (y(s),\psi(-s)),
\nonumber
\end{eqnarray}
or equivalently 
\begin{eqnarray}
z'(t)+K(t)z(t)=x'(t),\;\;z(0)=x(0), 
\nonumber
\end{eqnarray}
which can be considered as a linear inhomogeneous ODE in $z$. For such ODE the variation
of constants formula yields the $C^1[0,h]$-function
\begin{eqnarray}
z(t)=x(0)e^{-\int_0^tK(s)ds}+\int_0^te^{-\int_s^tK(\s)d\s}x'(s)ds.
\nonumber
\end{eqnarray}
One can use this expression to show bijectivity of $\d_2G(\psi,y)$ and as a definition of the
inverse. \qed\bigskip

Now, we fix $\ops\in M$ and denote by $\oy(\ops)$ the solution of (\ref{eq44}). By Lemma's
 \ref{lem6} and \ref{lem7} and the implicit function theorem we know that (\ref{eq44}) has
 solutions $y(\psi)$ with $(y,\psi)$ in a neighborhood of $(\oy,\ops)$. These solutions coincide with
 the solutions obtained via the Picard-Lindel\"of theorem. The new conclusion is that 
 $\psi\longmapsto y(\psi)$ is differentiable in a neighborhood of $\ops$. In this neighborhood we can
 differentiate the equation $G(\psi,y(\psi))=x_2$ with respect
to $\psi$ and apply it to $\chi$ from which we conclude that $Dy(\psi):C^1[-h,0]\_>C^1[0,h]$,
\begin{eqnarray}
&&Dy(\psi)\chi(t)=-\d_2 G (\psi, y(\psi))^{-1}\d_1G (\psi, y(\psi))\chi(t)
\nonumber\\
&=&-\int_0^te^{-\int_s^t\d_1 g (y(\s,\psi),\psi(-\s))d\s}\d_2 g (y(s,\psi),\psi(-s))\chi(-s)ds.
\label{eq19}
\end{eqnarray}
This formula can also define an extension, such that we can summarize as
\begin{lemma}
The map $y:M\_>C^1[0,h]$ satisfies (S1-S2), i.e., $y$ is $C^1$ with $Dy$ and $D_ey$
given by the right hand side of (\ref{eq19}).
\end{lemma}
To show that $y$ satisfies (S3), we prove
\begin{lemma}
Let $k:J\times I\_>\R$, $l:\R\_>\R$ and $m:C^1([0,h],J)\times M\_>C[0,h]$ be 
continuous maps. Then the following maps are continuous. 
\begin{eqnarray}
C^1([0,h],J)\times M&\_>&C[0,h],
\nonumber\\
(z,\psi)&\longmapsto&[t\longmapsto k(z(t),\psi(-t))],
\label{eq54}\\
C^1([0,h],J)\times M\times C[-h,0]&\_>&C[0,h],
\nonumber\\
(z,\psi,\chi)&\longmapsto&[t\longmapsto \int_0^{t}m(z,\psi)(s)\chi(-s)ds],
\nonumber
\\
\label{eq53}\\
C^1([0,h],J)\times M&\_>&C[0,h],
\nonumber\\
(z,\psi)&\longmapsto&[t\longmapsto l(\int_0^{t}k(z(s),\psi(-s))ds)].
\nonumber\\
\label{eq55}
\end{eqnarray}
\end{lemma}
\Proof Continuity of the map (\ref{eq54}) follows from Lemma 
\ref{lem12}. Next, 
\begin{eqnarray}
&&|\int_0^t m(z,\psi)(s)\chi(-s)ds-\int_0^t m(\oz,\ops)(s)\och(-s)ds|
\nonumber\\
&\le&h\|\chi\|\|m(z,\psi)-m(\oz,\ops)\|+\int_0^h |m(\oz,\ops)(s)|ds\|\chi-\och\|,
\nonumber
\end{eqnarray}
which implies continuity of (\ref{eq53}).
Let now $(\oz,\ops), (z,\psi)\in C^1([0,h],J)\times M$. 
By continuity of $t\mapsto k(z(t),\psi(-t))$, we get for $t\in [0,h]$
\begin{eqnarray}
|\int_0^tk(z(s),\psi(-s))ds|\le h\sup_{s\in[0,h]}|k(z(s),\psi(-s))|.
\nonumber
\end{eqnarray}
Then, by continuity of (\ref{eq54}) there exist $\de=\de(\oz,\ops)>0$, $L=L(\oz,\ops)>0$
such that
\begin{eqnarray}
|\int_0^t k(z(s),\psi(-s))ds|\le L,{\rm\;for\;all}\; t\in[0,h], (z,\psi)\in 
B((\oz,\ops),\de). 
\nonumber
\end{eqnarray}
Next, $l|_{[-L,L]}$ is uniformly continuous. To show continuity of (\ref{eq55}) it is thus
sufficient  to show continuity of 
\begin{eqnarray}
\begin{array}{ccc}
C^1([0,h],J)\times M&\_>&C[0,h],
\\
(z,\psi)&\longmapsto&[t\longmapsto\int_0^{t}k(z(s),\psi(-s))ds].
\end{array}
\nonumber
\end{eqnarray}
Continuity of (\ref{eq53}) implies this continuity with $m(z,\psi)(s):=k(z(s),\psi(-s))$
and $\chi\equiv1$. 
\qed\bigskip

Continuity of (\ref{eq54}--\ref{eq55}) helps to show
\begin{lemma}\label{lem15}
The map 
\begin{eqnarray}
&&
\begin{array}{ccc}
M\times C[-h,0]&\_>&C[0,h],
\\
(\psi,\chi)&\longmapsto& D_ey(\psi)\chi,
\end{array}
\nonumber\\
&&D_ey(\psi)\chi(t)=
-\int_0^{t}e^{-\int_s^t \d_1 g (y(\s,\psi),\psi(-\s))d\s}\d_2 g (y(s,\psi),\psi(-s))\chi(-s)ds
\nonumber\\
\label{eq57}
\end{eqnarray}
is continuous. 
\end{lemma}
\Proof We first show that
\begin{eqnarray}
\begin{array}{ccc}
C^1([0,h],J)\times M\times C[-h,0]&\_>&C[0,h],
\\
(z,\psi,\chi)&\longmapsto&[t\longmapsto
\int_0^{t}e^{-\int_s^t \d_1 g (z(\s),\psi(-\s))d\s}
\\
&&\cdot \d_2 g (z(s),\psi(-s))\chi(-s)ds]
\end{array}
\label{eq58}
\end{eqnarray}
is continuous. Then (\ref{eq57}) is continuous as a composition of the above map with
the continuous map
\begin{eqnarray}
\begin{array}{ccc}
M\times C[-h,0]&\_>&C^1([0,h],J)\times M\times C[-h,0],
\\
(\psi,\chi)&\longmapsto&
(y(\psi),\psi,\chi).
\end{array}
\nonumber
\end{eqnarray}
Continuity of (\ref{eq58}) follows as continuity of a product if we show continuity of the maps
\begin{eqnarray}
&&C^1([0,h],J)\times M\times C[-h,0]\_>C[0,h],
\nonumber\\
&&(z,\psi,\chi)\longmapsto
[t\longmapsto e^{-\int_0^t \d_1 g (z(\s),\psi(-\s))d\s}]
\label{eq59}\\
&&(z,\psi,\chi)\longmapsto[t\longmapsto
\int_0^t e^{\int_0^s \d_1 g (z(\s),\psi(-\s))d\s}\d_2 g (z(s),\psi(-s))\chi(-s)ds].
\nonumber\\
\label{eq60}
\end{eqnarray}
Continuity of (\ref{eq59}) and continuity of
\begin{eqnarray}
C^1([0,h],J)\times M&\_>&C[0,h],
\nonumber\\
(z,\psi)&\longmapsto&
[t\longmapsto e^{\int_0^t \d_1 g (z(\s),\psi(-\s))d\s}]
\label{eq61}
\end{eqnarray}
follow from the continuity of (\ref{eq55}). Continuity of (\ref{eq61})
together with continuity of (\ref{eq54}) imply continuity of 
\begin{eqnarray}
(z,\psi)&\longmapsto&
[t\longmapsto e^{\int_0^t \d_1 g (z(\s),\psi(-\s))d\s}\d_2 g (z(t),\psi(-t))].
\nonumber
\end{eqnarray}
This, together with continuity of (\ref{eq53}), implies continuity of (\ref{eq60}). 
\qed\bigskip

Now we can come back to our task of computing $D\t(\psi)$ and rewrite (\ref{eq9}) as 
\begin{eqnarray}
x_2-x_1=X(\psi,\t),
\label{eq45}
\end{eqnarray}
introducing
\begin{eqnarray}
X:M\times(0,h) &\longrightarrow&\R
\nonumber\\
(\psi,\t)&\longmapsto&X(\psi,\t):=\int_0^\t g(y(s,\psi),\psi(-s))ds,
\label{eq20}
\end{eqnarray}
where $y(s,\psi)$ is obtained from the Picard-Lindel\"of theorem. We show via the chain 
rule that $X$ is $C^1$. Let $p$ denote the projection of a vector on its first component and
$id_E$ the identity map on a set $E$. We then consider the map
\begin{eqnarray}
&&(y\circ p,id_{C^1\times \R}):M\times(0,h) \longrightarrow C^1([0,h],J)\times M\times
(0,h),
\nonumber\\
&&(y\circ p,id_{C^1\times \R})(\psi,\t)=(y(\psi),\psi,\t).
\nonumber
\end{eqnarray}
As $y$ is $C^1$, this map is also $C^1$ with derivative
\begin{eqnarray}
D(y\circ p,id_{C^1\times \R}):C^1\times \R\_>C^1([0,h],\R)\times C^1\times \R,
\nonumber\\
D(y\circ p,id_{C^1\times \R})(\psi,\t)(\chi,s)=(Dy(\psi)\chi,\chi,s).
\nonumber
\end{eqnarray}
Next, for continuous $k:J\times I\_>\R$, we define an integral operator
\begin{eqnarray}
I^k:C^1([0,h],J)\times M\times(0,h)& \longrightarrow&\R,
\nonumber\\
(z,\psi,\t)&\longmapsto&\int_0^\t k(z(s),\psi(-s))ds.
\label{eq74}
\end{eqnarray}
Then we can decompose $X$ in the way that 
\begin{eqnarray}
I^g\circ(y\circ p,id_{C^1\times\R})(\psi,\t)&=&I^g(y(\psi),\psi,\t)=\int_0^\t g(y(s,\psi),\psi(-s))ds
\nonumber\\
&=&X(\psi,\t).
\label{eq6}
\end{eqnarray}
We denote by ${\cal L}(E,F)$ the Banach space of bounded linear operators
from a Banach space $E$ into a Banach space $F$. 
To show that $I^g$ is $C^1$ and for later use, we show
\begin{lemma}\label{lem13}
Let $k:J\times I\_>\R$ be continuous. Then the following operators are continuous:
\begin{eqnarray}
\begin{array}{ccc}
C^1([0,h],J)\times M\times(0,h)& \longrightarrow&{\cal L}(C[0,h],\R),
\\
(z,\psi,\t)&\longmapsto&L^k(z,\psi,\t),
\end{array}
\label{eq48}\\
\begin{array}{ccc}
C^1([0,h],J)\times M\times(0,h)& \longrightarrow&{\cal L}(C^1[0,h],\R),
\\
(z,\psi,\t)&\longmapsto&L^k(z,\psi,\t),
\end{array}\label{eq49}\\
\begin{array}{ccc}
C^1([0,h],J)\times M\times(0,h)\times C[0,h]& \longrightarrow&\R,
\\
(z,\psi,\t,\chi)&\longmapsto&L^k(z,\psi,\t)\chi,
\end{array}
\label{eq65}
\\
\begin{array}{ccc}
C^1([0,h],J)\times M\times(0,h)& \longrightarrow&{\cal L}(C[-h,0],\R),
\\
(z,\psi,\t)&\longmapsto&J^k(z,\psi,\t),
\end{array}\label{eq66}
\\
\begin{array}{ccc}
C^1([0,h],J)\times M\times(0,h)& \longrightarrow&{\cal L}(C^1[-h,0],\R),
\\
(z,\psi,\t)&\longmapsto&J^k(z,\psi,\t),
\end{array}
\label{eq50}
\\
\begin{array}{ccc}
C^1([0,h],J)\times M\times(0,h)\times C[-h,0]& \longrightarrow&\R,
\\
(z,\psi,\t,\chi)&\longmapsto&J^k(z,\psi,\t)\chi,
\end{array}
\label{eq51}
\end{eqnarray}
where  $L^k(z,\psi,\t)\chi:=\int_0^\t k(z(s),\psi(-s))\chi(s)ds$ \\
and 
$J^k(z,\psi,\t)\chi:=\int_0^\t k(z(s),\psi(-s))\chi(-s)ds$.
\end{lemma}
\Proof We first show continuity of the operator (\ref{eq48}). One has
\begin{eqnarray}
&&|[L^k(z,\psi,\t)-L^k(\oz,\ops,\ot)]\chi|
\nonumber\\
&\le&
\|\chi\|[|\int_\ot^\t|k(z(s),\psi(-s))ds|
\nonumber\\
&&+\int_0^h|k(z(s),\psi(-s))-k(\oz(s),\ops(-s))|ds].
\nonumber
\end{eqnarray}
By continuity of $k$ and the compactness result in Lemma \ref{lem12} we know that 
$k(z(s),\psi(-s))$ is bounded for $(z,\psi)$ in a neighborhood of $(\oz,\ops)$.
It follows that in this neighborhood for some $K_1\ge 0$
\begin{eqnarray}
&&|[L^k(z,\psi,\t)-L^k(\oz,\ops,\ot)]\chi|
\nonumber\\
&\le&\|\chi\|[K_1|\t-\ot|+\int_0^h|k(z(s),\psi(-s))-k(\oz(s),\ops(-s))|ds].
\nonumber
\end{eqnarray}
The integral tends to zero if $\|z-\oz\|^1$ and $\|\psi-\ops\|^1$ tend to zero
by Lemma \ref{lem12}. Hence
\begin{eqnarray}
\|L^k(z,\psi,\t)-L^k(\oz,\ops,\ot)\|_{op}
\nonumber
\end{eqnarray}
tends to zero if $|\t-\ot|$, $\|z-\oz\|^1$ and $\|\psi-\ops\|^1$ tend to zero. Thus (\ref{eq48}) is continuous. We can conclude 
immediately that the operators defined in (\ref{eq49}--\ref{eq65}) are continuous. 
The remaining continuity statements can be shown analogously.
\qed\bigskip

\begin{lemma}\label{lem10}
Let $k:J\times I\_>\R$ be $C^1$. Then
$I^k$ defined by (\ref{eq74}) is $C^1$.  For $(z,\psi,\t)\in C^1([0,h],J)\times M\times(0,h)$ the derivative is given by
\begin{eqnarray}
&&DI^k(z,\psi,\t)=
(L^{\d_1k}(z,\psi,\t),J^{\d_2k}(z,\psi,\t),k(z(\t),\psi(-\t)))
\nonumber\\
&&
:C^1[0,h]\times C^1[-h,0]\times\R\_>\R,
\end{eqnarray}
 where $L^{\d_1k}$ and $J^{\d_2 k}$ are as in Lemma \ref{lem13} with $k$ replaced by the respective 
partial derivatives.
\end{lemma}
\Proof  First, 
\begin{eqnarray}
&&|I^k(z,\psi,\t)-I^k(\oz,\psi,\t)-L^{\d_1k}(\oz,\psi,\t)(z-\oz)|
\nonumber\\
&\le&\int_0^h |k(z(s),\psi(-s))-k(\oz(s),\psi(-s))
\nonumber\\
&&-\d_1k(\oz(s),\psi(-s))(z(s)-\oz(s))|ds.
\nonumber
\end{eqnarray}
For fixed $\psi$ and $\t$ the right hand side is $o(\|z-\oz\|^1)$ as $\|z-\oz\|^1\rightarrow 0$ 
by Lemma \ref{lem11} (a). 
Next, continuity of $(z,\psi,\t)\longmapsto L^{\d_1k}(z,\psi,\t)$ follows by continuity of 
(\ref{eq49}). We have shown partial differentiability in the first argument and continuity of 
the partial derivative.  The corresponding result for the second argument is proven similarly, but
using Lemma \ref{lem11} (b) and continuity of (\ref{eq50}). The corresponding result for the third argument is clear. The statement follows.
\qed\bigskip

We can apply the chain rule in (\ref{eq6}) and deduce
\begin{lemma}\label{lem16}
The functional $X$ defined in (\ref{eq20}) is $C^1$ with derivative $DX=(\d_1 X,\d_2 X)$ where 
\begin{eqnarray}
&&\d_1 X(\psi,\t):C^1[-h,0]\_>\R,
\nonumber\\
&&
\d_1 X(\psi,\t)\chi
=L^{\d_1g}(y(\psi),\psi,\t)Dy(\psi)\chi+J^{\d_2g}(y(\psi),\psi,\t)\chi,
\nonumber
\\
&&\d_2 X(\psi,\t):\R\_>\R;\;\d_2 X(\psi,\t)1=g(y(\t,\psi),\psi(-\t)).
\nonumber
\end{eqnarray}
\end{lemma}

Now note that $\d_2 X(\psi,\t)1=g(y(\t,\psi),\psi(-\t))\neq 0$. 
Hence, by the implicit function theorem applied to (\ref{eq45}) the operator  $\psi\longmapsto\t(\psi)$ is $C^1$. We obtain the derivative
by differentiating $x_2-x_1=X(\psi,\t(\psi))$ with respect to $\psi$, which yields
\begin{eqnarray}
0&=&\frac{d}{d\psi}X(\psi,\t(\psi))\chi
\nonumber\\
&=&L^{\d_1g}(y(\psi),\psi,\t(\psi))Dy(\psi)\chi+J^{\d_2g}(y(\psi),\psi,\t(\psi))\chi
\nonumber\\
&&+g(y(\t(\psi),\psi),\psi(-\t(\psi)))D\t(\psi)\chi.
\nonumber
\end{eqnarray}
Hence with $y(\t(\psi),\psi)=x_1$, we get 
\begin{eqnarray}
&&D\t(\psi)\chi
\nonumber\\
&=&-\frac{L^{\d_1g}(y(\psi),\psi,\t(\psi))Dy(\psi)\chi+J^{\d_2g}(y(\psi),\psi,\t(\psi))\chi}
{g(x_1,\psi(-\t(\psi)))},
\nonumber\\
\label{eq29}
\end{eqnarray}
with $Dy$ as in (\ref{eq19}). Next, we define an extension $D_e\t$
of $D\t$ on $C[-h,0]$ by (\ref{eq29}) with $Dy$ replaced by $D_ey$. Then
we have obtained
\begin{lemma}\label{lem24}
The functional $\t$ satisfies (S1-S2), i.e., $\t:M\_>[0,h]$ is $C^1$ with $D\t$ defined by
(\ref{eq29}) and $D_e\t$ is defined by (\ref{eq29}) with $Dy$ replaced by $D_ey$. 
\end{lemma}
To show that (S3) holds, we prove the following results. 
\begin{lemma}\label{lem8}
Suppose that $k:J\times I\_>\R$ is continuous. Then the following 
functionals are continuous. 
\begin{eqnarray}
\begin{array}{ccc}
M\times C[0,h]&\_>&\R,
\\
(\psi,\chi)&\longmapsto&L^k(y(\psi),\psi,\t(\psi))\chi,
\end{array}
\label{eq67}\\
\begin{array}{ccc}
M\times C[-h,0]&\_>&\R,
\\
(\psi,\chi)&\longmapsto&J^k(y(\psi),\psi,\t(\psi))\chi.
\end{array}
\label{eq52}
\end{eqnarray}
\end{lemma}
\Proof We denote the operator defined in (\ref{eq65}) by $K^k$, i.e., 
\begin{eqnarray}
&&K^k:C^1([0,h],J)\times M\times(0,h)\times C[0,h]\_>\R,
\nonumber\\
&&K^k(z,\psi,\t,\chi):=L^k(z,\psi,\t)\chi.
\nonumber
\end{eqnarray}
Next we denote by $q$ the projection of a vector on its second component. Then
\begin{eqnarray}
((y,id_{C^1},\t)\circ p,q)(\psi,\chi)=(y(\psi),\psi,\t(\psi),\chi).
\nonumber
\end{eqnarray}
Hence, 
\begin{eqnarray}
&&K^k\circ ((y,id_{C^1},\t)\circ p,q)(\psi,\chi)=K^k(y(\psi),\psi,\t(\psi),\chi)
\nonumber\\
&=&L^k(y(\psi),\psi,\t(\psi))\chi.
\nonumber
\end{eqnarray}
Continuity of $K^k$ is continuity of (\ref{eq65}) and $((y,id_{C^1},\t)\circ p,q)$ is continuous by our earlier results. Hence 
(\ref{eq67}) is continuous as a composition of continuous maps. Continuity of (\ref{eq52}) follows analogously if one uses continuity of (\ref{eq51}). 
\qed\bigskip

\begin{lemma}\label{lem14}
Let $k:J\times I\_>\R$ be continuous. Then the following functional is continuous. 
\begin{eqnarray}
\begin{array}{ccc}
M\times C[-h,0]&\_>&\R,
\\
(\psi,\chi)&\longmapsto&L^k(y(\psi),\psi,\t(\psi))D_ey(\psi)\chi.
\end{array}
\label{eq46}
\end{eqnarray}
\end{lemma}
\Proof 
We define
\begin{eqnarray}
&F_1:M\times C[0,h]\_>\R;&
F_1(\psi,\chi):=L^k(y(\psi),\psi,\t(\psi))\chi,
\nonumber\\
&F_2:M\times C[-h,0]\_>M\times C[0,h];&
F_2(\psi,\chi):=(\psi, D_ey(\psi)\chi),
\nonumber
\end{eqnarray}
such that (\ref{eq46}) is $F_1\circ F_2$. 
Continuity of $F_1$ is continuity of  (\ref{eq67}) and $F_2$ is continuous by Lemma 
\ref{lem15}, thus the statement follows. 
\qed\bigskip

Then we can prove
\begin{lemma}\label{lem17}
The functional $\t:M\_>[0,h]$ satisfies (S3), i.e., 
\begin{eqnarray}
&M\times C[-h,0]\longrightarrow\R
;&(\psi,\chi)\longmapsto D_e\t(\psi)\chi,
\nonumber
\end{eqnarray}
with $D_e\t$ defined in Lemma \ref{lem24} is continuous.
\end{lemma}
\Proof 
Continuity of
\begin{eqnarray}
M\times C[-h,0]&\longrightarrow&\R
\nonumber\\
(\psi,\chi)&\longmapsto&J^{\d_2g}(y(\psi),\psi,\t(\psi))\chi
\nonumber\\
(\psi,\chi)&\longmapsto&L^{\d_1g}(y(\psi),\psi,\t(\psi))D_ey(\psi)\chi
\nonumber
\end{eqnarray}
is implied by continuity of (\ref{eq52}) and (\ref{eq46}) respectively.  Also
\begin{eqnarray}
M\times C[-h,0]\longrightarrow\R
;
(\psi,\chi)&\longmapsto&\frac{1}{g(x_1,\psi(-\t(\psi)))}
\nonumber
\end{eqnarray}
is continuous. Hence, $(\psi,\chi)\longmapsto D_e\t(\psi)\chi$ is continuous. 
\qed\bigskip

We can now summarize the previous results as
\begin{theorem}\label{theo8}
Suppose that $g$ satisfies (G). Then the functional $\t$ describing 
the delay satisfies the smoothness property (S). 
\end{theorem}
We will now show that also $\cF$  satisfies (S). 
\begin{theorem}\label{theo6}
Suppose that $g$ satisfies (G) and that moreover $d:J\times I\_>\R$ is $C^1$. Then 
\begin{eqnarray}
\cF(\psi)=e^{\int_0^{\t(\psi)} d(y(s,\psi),\psi(-s))ds}=exp\circ I^d\circ(y,id_M,\t)(\psi)
\nonumber
\end{eqnarray}
fulfills (S). In particular 
\begin{eqnarray}
D\cF(\psi)&=&e^{I^d(y(\psi),\psi,\t(\psi))}[L^{\d_1d}(y(\psi),\psi,\t(\psi))Dy(\psi)+
J^{\d_2d}(y(\psi),\psi,\t(\psi))
\nonumber\\
&&+d(x_1,\psi(-\t(\psi)))D\t(\psi)]
\label{eq33}
\end{eqnarray}
with extension 
\[
D_e\cF(\psi):C[-h,0]\_>\R
\]
 of $D\cF(\psi)$ is defined with the right hand side of (\ref{eq33}) replacing $D\t(\psi)$ and $Dy(\psi)$
with $D_e\t(\psi)$ and $D_ey(\psi)$. 
\end{theorem}
\Proof First, $I^d$ is $C^1$ by Lemma \ref{lem10} applied to $d$. Hence, $\cF$ is $C^1$ by the chain rule. We should show that 
\begin{eqnarray}
M\times C[-h,0]\longrightarrow\R;
(\psi,\chi)\longmapsto D_e\cF(\psi)\chi
\nonumber
\end{eqnarray}
is continuous. First, we know that the map
\begin{eqnarray}
M\times C[-h,0]\longrightarrow\R;(\psi,\chi)\longmapsto e^{I^d(y(\psi),\psi,\t(\psi))}
\nonumber
\end{eqnarray}
is continuous as a composition by continuity of $y$, $\t$, $exp$ and $I^d$. Next, 
\begin{eqnarray}
(\psi,\chi)&\longmapsto&d(x_1,\psi(-\t(\psi)))D_e\t(\psi)\chi
\nonumber
\end{eqnarray}
is continuous by Lemma \ref{lem17}. Moreover continuity of
\begin{eqnarray}
(\psi,\chi)&\longmapsto&L^{\d_1d}(y(\psi),\psi,\t(\psi))D_ey(\psi)\chi,
\label{eq36}\\
(\psi,\chi)&\longmapsto&J^{\d_2d}(y(\psi),\psi,\t(\psi))\chi
\label{eq37}
\end{eqnarray}
is implied by continuity of (\ref{eq46}) applied to $\d_1d$ and continuity of (\ref{eq52}) 
applied to $\d_2d$ respectively.  In summary we have shown
that $\cF$ satisfies (S). 
\qed\bigskip

%
\subsubsection{Global existence}
We should guarantee that $\t$, as defined in the Picard-Lindel\"of theorem, and 
$\cF(\psi)=e^{-\int_0^{\t(\psi)}d(y(s,\psi),\psi(-s))ds}$ fulfill the hypotheses for global 
existence. In 
the following we will see how the boundedness properties that have been assumed for the
derivatives in Theorem \ref{theo5} (d) will be used. For the solution $y$, that we 
have obtained via the Picard-Lindel\"of theorem, we can show
\begin{lemma}\label{lem2}
Suppose that $g$ satisfies (G), $k:J\times I\_>\R$ is $C^1$ and that
$\d_2g(\oB\times A)$ and $\d_ik(\oB\times A)$, $i=1,2$, are bounded for any bounded
$A\subset I$. Then for any $C$-bounded $B\subset M$ there exist some $K_B,L_B, M_B\ge 0$ such that
\begin{eqnarray}
\int_0^h|k(y(s,\psi),\chi(-s))-k(y(s,\psi),\psi(-s))|ds&\le& M_B\|\psi-\chi\|,
\nonumber\\
\int_0^h|y(s,\chi)-y(s,\psi)|ds&\le& K_B\|\chi-\psi\|,
\nonumber\\
\int_0^h|k(y(s,\psi),\psi(-s))-k(y(s,\chi),\chi(-s))|ds&\le& L_B\|\psi-\chi\|.
\label{eq76}
\end{eqnarray}
for all $\psi, \chi\in B$. 
\end{lemma}
\Proof 
For some $\th_s\in[0,1]$
\begin{eqnarray}
&&|k(y(s,\psi),\chi(-s))-k(y(s,\psi),\psi(-s))|
\nonumber\\
&=&|\d_2k(y(s,\psi),\chi(-s)+\th_s(\psi(-s)-\chi(-s)))||\chi(-s)-\psi(-s)|.
\nonumber
\end{eqnarray}
Then the boundedness property of $\d_2k$ implies the first statement. For the bounded set 
$A:=\{\chi(-s):\;\chi\in B,s\in[0,h]\}$ 
we  define 
\[
K_2:=\sup_{(x,y)\in\oB\times A}|\d_1k(x,y)|
\]
and get
\begin{eqnarray}
|k(y(s,\chi),\chi(-s))-k(y(s,\psi),\chi(-s))|
\le K_2|y(s,\chi)-y(s,\psi)|. 
\label{eq63}
\end{eqnarray}
Then, by (\ref{eq10}), the previous estimate and 
the first statement applied to $g$ instead of $k$
\begin{eqnarray}
&&\int_0^h|y(s,\chi)-y(s,\psi)|ds
\nonumber\\
&=&\int_0^h|\int_0^sg(y(\s,\chi),\chi(-\s))-g(y(\s,\psi),\psi(-\s))d\s|ds
\nonumber\\
&\le&h\int_0^h|g(y(s,\chi),\chi(-s))-g(y(s,\psi),\psi(-s))|ds
\nonumber\\
&\le&h\int_0^h|g(y(s,\chi),\chi(-s))-g(y(s,\psi),\chi(-s))|ds
\nonumber\\
&&+h\int_0^h|g(y(s,\psi),\chi(-s))-g(y(s,\psi),\psi(-s))|ds
\nonumber\\
&\le&hK_2\int_0^h|y(s,\chi)-y(s,\psi)|ds
+hM_B\|\chi-\psi\|,
\nonumber
\end{eqnarray}
where $M_B\ge 0$ is some constant
that exists by the first statement. Hence, as $h=\frac{b}{K}$, we have
\[
hK_2=\frac{b}{K}\sup_{(x,y)\in\oB\times A}|\d_1 g (x,y)|<1
\]
by (G2) and get that
\begin{eqnarray}
\int_0^h|y(s,\chi)-y(s,\psi)|ds
\le \frac{hM_B}{1-hK_2}\|\chi-\psi\|,
\nonumber
\end{eqnarray}
which proves the second statement. 
If we combine this result with (\ref{eq63}),  we get that 
\begin{eqnarray}
&&\int_0^h|k(y(s,\psi),\chi(-s))-k(y(s,\chi),\chi(-s))|ds
\le K_3\|\chi-\psi\|
\nonumber
\end{eqnarray}
for some $K_3\ge 0$. The third statement follows from the first statement
and the previous estimate. 
\qed\bigskip

Without further assumptions on $g$, we get
\begin{lemma}\label{lem9}
Suppose that $g$ satisfies (G) and that $\d_2g(\oB\times A)$ is
bounded, whenever $A\subset I$ is bounded. Then for any $C$-bounded $B\subset M$, there exists some $L_B\ge 0$, such that 
\[
|\t(\psi_1)-\t(\psi_2)|\le L_B\|\psi_1-\psi_2\|, \;{\rm for}\;{\rm all}\;\psi_1,\psi_2\in B.
\]
In particular, $\t$ is (Lb). 
\end{lemma}

\Proof Let $B\subset M$ be $C$-bounded, $\psi,\chi\in B$. By (\ref{eq9}) we get that
\begin{eqnarray}
\int_0^{\t(\psi)}g(y(s,\psi),\psi(-s))ds=\int_0^{\t(\chi)}g(y(s,\chi),\chi(-s))ds
\nonumber
\end{eqnarray}
or, equivalently,
\begin{eqnarray}
&&\int^{\t(\psi)}_{\t(\chi)}g(y(s,\psi),\psi(-s))ds
\nonumber\\
&=&
\int_0^{\t(\chi)}g(y(s,\chi),\chi(-s))-g(y(s,\psi),\psi(-s))ds.
\nonumber
\end{eqnarray}
This implies that by (G3)
\begin{eqnarray}
&&|\t(\psi)-\t(\chi)|\le\frac{1}{\e}
\int_0^{h}|g(y(s,\chi),\chi(-s))-g(y(s,\psi),\psi(-s))|ds
\nonumber\\
&\le& K\|\psi-\chi\|
\nonumber
\end{eqnarray}
for some $K\ge 0$ by (\ref{eq76}) applied to $g$ in place of $k$, which proves the statement. 
\qed\bigskip

Next, note that the boundedness properties of $\t$ and $d$ shown and assumed in 
Theorem \ref{theo5} (a) and (d) respectively, imply the boundedness of $\cF$. Moreover, we can prove
\begin{lemma}\label{lem25}
Suppose that $g$ satisfies (G), $d:J\times I\longrightarrow \R$ is $C^1$ and the sets
\[
d(\oB\times I),\;\d_2 g(\oB\times A), \;\d_i d(\oB\times A), \;i=1,2
\]
are bounded, whenever $A\subset I$ is bounded. Then for any $C$-bounded $B\subset M$, there exists some $L_B\ge 0$, such that 
\[
|\cF(\psi_1)-\cF(\psi_2)|\le L_B\|\psi_1-\psi_2\|, \;{\rm for}\;{\rm all}\;\psi_1,\psi_2\in B.
\]
In particular,  $\cF$ is (Lb). 
\end{lemma}
\Proof Let $B\subset M$ be $C$-bounded, $\psi,\chi\in B$. Then for some $\th_s\in[0,1]$
\begin{eqnarray}
|\cF(\psi)-\cF(\chi)|&=&|e^{I^d(y(\psi),\psi,\t(\psi))}-e^{I^d(y(\chi),\chi,\t(\chi))}|
\nonumber\\
&=&
e^{I^d(y(\psi),\psi,\t(\psi))+\th_s[I^d(y(\chi),\chi,\t(\chi))-I^d(y(\psi),\psi,\t(\psi)]}
\nonumber\\
&&\cdot
|I^d(y(\psi),\psi,\t(\psi))-I^d(y(\chi),\chi,\t(\chi))|.
\nonumber
\end{eqnarray}
Next, note that the boundedness condition for $d$ ensures that 
\begin{eqnarray}
B\_>\R,\;\psi\longmapsto I^d(y(\psi),\psi,\t(\psi))
\label{eq64}
\end{eqnarray}
is bounded. It is sufficient to show that in addition this map satisfies the discussed Lipschitz property. 
We have
\begin{eqnarray}
&&|I^d(y(\psi),\psi,\t(\psi))-I^d(y(\chi),\chi,\t(\chi))|
\nonumber\\
&=&|\int_0^{\t(\psi)}d(y(s,\psi),\psi(-s))ds-\int_0^{\t(\chi)}d(y(s,\chi),\chi
(-s))ds|
\nonumber\\
&\le&|\int_{\t(\chi)}^{\t(\psi)}d(y(s,\psi),\psi(-s))ds|
\nonumber\\
&&
+\int_0^{\t(\chi)}|d(y(s,\psi),\psi(-s))-d(y(s,\chi),\chi(-s))|ds.
\nonumber
\end{eqnarray}
The first integral is dominated by
\begin{eqnarray}
 K_1|\t(\chi)-\t(\psi)|\le K_2\|\psi-\chi\|
\nonumber
\end{eqnarray}
for some $K_i\ge 0$, $i=1,2$, by the boundedness property of $d$ and Lemma 
\ref{lem9}. The second integral is dominated by
\begin{eqnarray}
K_3 \|\psi-\chi\|
\nonumber
\end{eqnarray}
for some $K_3\ge 0$ by (\ref{eq76}) applied to $d$ in place of $k$.
In summary, (\ref{eq64})
satisfies the Lipschitz property, hence so does $\cF$. 
\qed\bigskip

\bigskip
\noindent{\bf Proof of Theorem \ref{theo5}.} 
First, (G1-2) imply conditions (i-ii) of the Picard-Lindel\"of theorem, which thus can
be applied to guarantee existence and uniqueness of $y$ and $\t$ as stated. Next, by 
Theorem \ref{theo8}, $\t$ satisfies (S), thus (a) follows. By Theorem \ref{theo6}, $\cF$ also satisfies (S), hence (b) follows. Then (c) follows by Theorem \ref{theo4}. 
Next, by Lemmas \ref{lem9} and \ref{lem25}, $\t$ and $\cF$ are (Lb). Since we have 
guaranteed that $\cF$ is bounded, also the remaining conditions of Theorem \ref{theo4} (b) are satisfied, the theorem implies
statement (d) of Theorem \ref{theo5}.
\qed\bigskip

\section{Discussion}
In this paper we have analyzed global existence and uniqueness for a differential equation with implicitly defined delay with state dependence. 
We have elaborated a new sufficient criterion for global existence (Theorem \ref{theo1}
and Corollary \ref{cor1}) for differential equations with state dependent delay. The equation that we have analyzed
describes a model for the maturation process of a stem cell population. 
We have elaborated conditions, with the aim of keeping them minimal, for basic model ingredients, i.e., vital rates, that guarantee local and - via above mentioned criterion - global existence and uniqueness of solutions. 
It has become clear that a merely implicit definition of the delay complicates the verification of smoothness conditions (conditions (S) and (Lb)).  

For a function $f$ with real arguments and values and an $\R$-valued functional $r$ the state dependent delay differential equation 
\begin{eqnarray}
x'(t)=f(x(t-r(x_t)))
\label{eq77}
\end{eqnarray}
is thoroughly analyzed in \cite{Walther,Louihi,Mallet,Walther3}. In particular these papers
contain well-posedness and stability results. Stability of periodic solutions of the equation
\begin{eqnarray}
x'(t)=f(x(t),x(t-r(x(t)))),
\label{eq79}
\end{eqnarray}
is analyzed in \cite{Mallet2, Mallet3}. 
In \cite{Crauste} Adimy and coauthors present a detailed stability analysis for a model that, as ours, describes the maturation of stem cells. The resulting equation is (\ref{eq79}). The authors consider either general or explicitly given state dependencies of the delay. 

In \cite{Hu}  and \cite {Balanov} the authors analyze the existence of Hopf bifurcations for equations that could, for comparison, be written as
\begin{eqnarray}
x'(t)=f(x(t), x(t-r(x_t,t))) \;{\rm and}\; x'(t)=f(x(t), x(t-r(x_t)), r(x_t)),
\nonumber\\
\label{eq80}
\end{eqnarray}
respectively. The functionals $r$ that specify the delay in terms of the state are defined implicitly, but not via an ODE, as in our case. 

For comparison we may write the central equation of our studies, i.e., (\ref{eq68}--\ref{eq11}), as
\begin{eqnarray}
x'(t)=f(x_t,x(t-r(x_t)),r(x_t))
\label{eq78}
\end{eqnarray}
(where the direct $r$-dependence in the third argument comes in via the $\t$-dependence
of $\cF$ in (\ref{eq81})).
We first remark that other than (\ref{eq77}) and (\ref{eq79}) our equation is $\mathbb R^2$-
valued. Recall that the state of the system at time $t$ is $x_t$ and 
thus $x(t)$ is the state at time $t$ evaluated in zero. 
Then, other than in the equations in (\ref{eq77}--\ref{eq80}), we have in both  the first argument of $f$ as well as in the functional $r$ that describes the delay a general state-dependence rather than (some) state-dependencies evaluated in zero as in those
equations. 
Our equation (\ref{eq78}) can thus not be written in either of the forms given in
(\ref{eq77}--\ref{eq80}). On the other hand a one dimensional variant of (\ref{eq78}) is a generalization of (\ref{eq77}), (\ref{eq79}) and the second equation in (\ref{eq80}). 

In \cite{Alarcon} Alarc{\'o}n et al. considered ingredients of our model specified in terms of parameters. For this parameterization they have analyzed the possibility of a unique positive equilibrium and computed representations of it.
It is shown in \cite{Getto2, Nakata} that for ODE variants of our model there is the possibility of destabilization of equilibria via Hopf bifurcation and the emergence of oscillations. 
 As a future
project we plan the stability analysis of equilibria in the general setting of the present paper. 
One of the first goals is thus a linearization of our equation. For this purpose
we plan to investigate applicability of the theory for differentiable semiflows developed in \cite{Walther4}. 
Moreover, we would like to investigate the existence of periodic solutions and how they relate to biological mechanisms at the cell level. We hope that the theory for equations (\ref{eq77}--\ref{eq80}) developed in above mentioned references can be extended to our equation (\ref{eq78}) and that the top down approach to modeling ingredients that we have used here can be applied to establish also a stability analysis of our model.

Though for (\ref{eq77}--\ref{eq80}) global existence is established, 
in ways different from ours, our criterion for global existence may also be interesting
for further analysis of (applications of) those equations.

For estimations of trajectories we have used a variation of constants formula
for the deduction of which we have exploited a certain linearity of the population
equations, more precisely, the fact that they could be written in the form (\ref{eq69}). 
We remark that this form, based on mass action laws, arises naturally in population dynamical modeling. 

The model behind our equation incorporates physiological structure of individual
cells with respect to different stages of maturity. 
It is typical for structured population models that individual development, unless it is age, 
depends on the population state, which is the fact that leads to the implicitly defined delay with state-dependence here, see e.g. \cite{Claessen, Diekmann, Roos}. 
See \cite[Section 2.1]{Walther} for an example from physics (two-body problem)
involving implicitly defined state-dependent delays. 
A general class of structured population models was analyzed in \cite{Diekmann} (local existence) and \cite{Getto} (global existence and continuous dependence).
 Application of this approach to our model could start as follows: Consider an initial 
 population measure on the interval $[x_1,x_2]$. Define so called input functions such that 
 assuming that these are given, the future population measure on $[x_1,x_2]$ depends {\it 
 linearly} on the initial measure. With the operator mapping initial measures to future 
 measures one should then define and solve a fixed point problem for the input functions. 
 An intuitive approach would be to take our two state components $w$ and $v$ as input functions. The future 
 population measure could
be defined via integration along the characteristics. The output-input map, i.e., the fixed 
point map, however would then be evaluation of the population measure in $x_1$ and 
$x_2$, which cannot be defined via integration of a continuous function with respect to the 
measure as is required in \cite{Diekmann, Getto}. It hence is not clear how to define and 
solve the fixed point problem with this approach. 

Similarly, in general also a partial differential equation formulation leads to solving a fixed point problem, see however \cite{Doumic} for an analysis of a special case of the model described by (\ref{eq68}--\ref{eq11}) via a limiting argument for multi-compartment models.

\section*{Acknowledgements}
The research of Ph.G. is part of his project
``Delay equations and structured population dynamics'' funded by the DFG (Deutsche
Forschungsgemeinschaft). The author receives additional support by the Spanish Ministry of Economy and Competitiveness (MINECO) under project MTM 2010-18318.

Ph.G. would like to thank Hans-Otto Walther for valuable comments on several versions
of the manuscript as well as Eugen Stumpf for interesting discussions on this research. 





\end{document}